\documentclass[letterpaper,12pt,reqno]{amsart}
\usepackage{graphicx,amsmath,amssymb,amsthm}
\usepackage[left=1in,right=1i
n, top=1in, bottom=1in]{geometry}
\newtheorem{thm}{Theorem}[section]

\newtheorem{prop}[thm]{Proposition}
\newtheorem{thrm}[thm]{Theorem}
\newtheorem{lem}[thm]{Lemma}

\newtheorem{rmks}[thm]{Remarks}
\newtheorem{conj}[thm]{Conjecture}

\newcommand{\RR}{\mathbb{R}}
\newcommand{\QQ}{\mathbb{Q}}

\newcommand{\NN}{\mathbb{N}}
\newcommand{\ZZ}{\mathbb{Z}}
\newcommand{\FF}{\mathbb{F}}
\newcommand{\C}{{\mathcal{C}}}
\newcommand{\aut}{\rm Aut}
\newcommand{\FFp}{S,T\in\FF(P)}
\newcommand{\FFps}{S,T\in\FF(P)^*}
\newcommand{\FFpu}{\exists (u_1,\ldots, u_L)\in\FF(P)^*}
\newcommand{\wpab}{w(P,E_{a,b})}

\newcommand{\sumpchiprim}{\underset{\substack{p\leq X \\  p_i^-<p_{i+1}<p_{i}^+ \\ 1\leq i \leq L-1}}{{\sum }^*}}
\newcommand{\sumofallps}{\sum_{\substack{p\leq X \\ p_i^-<p_{i+1}<p_{i}^+ \\ 1\leq i \leq L-1}}}
\newcommand{\sumofallh}{\sum_{\substack{p_i^-<p_{i+1}<p_i^+ \\ 1\leq i \leq L-1}}\prod_{j=1}^L H(D(p_j,p_{j+1}))}

\newcommand{\sumofpfracsing}{\sum_{\substack{p\leq X \\ p_i^-<p_{i+1}<p_{i}^+ \\ 1\leq i \leq L-1}}\frac{1}{p_1\cdots p_L}}
\newcommand{\sumofpminusone}{\sum_{\substack{p\leq X \\ p_i^-<p_{i+1}<p_{i}^+ \\ 1\leq i \leq L-1}}\frac{1}{(p_1-1)\cdots (p_L-1)}}

\makeatletter
\def\imod#1{\allowbreak\mkern7mu({\operator@font mod}\,\,#1)}
\makeatother

\numberwithin{equation}{section}

\title[Amicable pairs and aliquot cycles on average]{Amicable pairs and aliquot cycles on average}
\author{James Parks}
\address{Department of Mathematics and Computer Science, Lethbridge University, 4401 University Drive, Lethbridge, AB, T1K 3M4, Canada}
\email{james.parks@uleth.ca}
\thanks{The author is supported by a PIMS Postdoctoral Fellowship.}

\date{\today}
\begin{document}

\begin{abstract}
Silverman and Stange defined the notion of an aliquot cycle of length $L$ for a fixed elliptic curve $E/\QQ$, 
and conjectured an order of magnitude for the function that counts such aliquot cycles. We show that the conjectured upper bound holds for 
the number of aliquot cycles on average over the family of all elliptic curves with short bounds on the size of the parameters in the family.
\end{abstract}

\maketitle

\section{Introduction}
Let $E$ be an elliptic curve defined over $\QQ$ and let $L \geq 2$ be a positive integer. For a prime $p$, let $a_p(E)$ denote the trace of the Frobenius automorphism. 
Silverman and Stange \cite{JSKS:1} defined an $L$-tuple $(p_1,\ldots,p_L)$ of distinct prime numbers to be an \textit{aliquot cycle} of length $L$ of $E$ 
if $E$ has good reduction at each prime $p_i$ and  $$\#{E}_{p_i}(\FF_{p_i})=p_i+1-a_{p_i}(E_{p_i})=p_{i+1} \quad {\rm for}\: 1\leq i\leq L,$$
where we set $p_{L+1}:=p_1$. Aliquot cycles of length $L=2$ are called \textit{amicable pairs}. These definitions can be interpreted as the 
elliptic curve analogues to the classically defined aliquot cycles. As observed in \cite[Remark 1.5]{JSKS:1} aliquot cycles arose naturally when Silverman and Stange
generalized Smyth's \cite{CS:1} results on index divisibility of Lucas sequences  to elliptic divisibility sequences. 

We are interested in the the distribution of aliquot cycles of a given length $L$ for a fixed elliptic curve $E/\QQ$. We define an aliquot cycle $(p_1,\ldots,p_L)$ to 
be \textit{normalized} if $p_1=\min\{p_i:1\leq i\leq L\}$. We consider the normalized aliquot cycle counting function 
\begin{equation*}
\pi_{E,L}(X):=\#\{(p_1,\ldots,p_L)\:{\rm is\:a\:normalized\:aliquot\:cycle}\:|\:p_1\leq X\}.
\end{equation*}
Silverman and Stange \cite{JSKS:1} used a heuristic argument to give the following conjecture for the behavior of $\pi_{E,L}(X)$.
\begin{conj}[\textbf{Silverman-Stange}] \label{silverstangeamconj}
Let $E/\QQ$ be an elliptic curve and let $L\geq 2$ be a positive integer. Assume that there are infinitely many primes $p_i$ such that 
$\#{E}_{p_i}(\FF_{p_i})$ is prime. Then as $X\rightarrow\infty$ we have that
\begin{align*}
\pi_{E,L}(X) &\asymp \frac{\sqrt{X}}{(\log X)^L} \quad {\rm if}\:E\:{\rm does\:not\:have\:complex\:multiplication\:(CM)},\\
\pi_{E,2}(X)&\sim  A_E\frac{X}{(\log X)^2} \quad {\rm if}\:E\:{\rm has\:CM},
\end{align*}
where the implied constants in $\asymp$ are both positive and depend only on $E$ and $L$ and $A_E$ is a precise positive constant.
\end{conj}

\begin{rmks}
{\em (i) We may interpret the case $L=1$ in Conjecture \ref{silverstangeamconj} as 
describing primes $p$ for which $\#{E}_{p}(\FF_{p})=p$. 
These primes are called \textit{anomalous primes} and were previously considered by Mazur \cite{BM:1}. In this case, 
Conjecture \ref{silverstangeamconj} is a special case of a conjecture of Lang and Trotter \cite{SLHT:1}.

(ii) Silverman and Stange \cite{JSKS:1} focused primarily on the CM case. They showed that if $E/\QQ$ has CM with $j$-invariant $j_E\neq 0$ 
then there are no normalized aliquot cycles of length $L\geq 3$ for primes $p\geq 5$. This implies that $\pi_{E,L}(X)=O(1)$. If $E$ has CM with $j_E=0$ then they showed that $E$ does not have any normalized aliquot triples $(p,q,r)$ with $p>7$. However, it is unknown if $\pi_{E,L}(X)=O(1)$ when $j_E=0$ and $L> 3$ and no conjecture is given in this case. Also, no formula is given for $A_E$ in Conjecture \ref{silverstangeamconj}.

(iii) We remark that for $1\leq i \leq L-1$, we have that 
\begin{equation} \label{pplusminus}
p_i^-:=p_i+1-2\sqrt{p_i}< p_{i+1}:=\# E_{p_i}(\FF_{p_i}) < p_i^+:=p_i+1+2\sqrt{p_i}
\end{equation}
by Hasse's Theorem (see \cite[Chapter V, Theorem 1.1]{JS:1}).}
\end{rmks}

Jones \cite{NJ:1} refined Conjecture \ref{silverstangeamconj} in the non-CM case. He gave a precise conjectural constant $C_{E,L}$ in the asymptotic formula for $\pi_{E,L}(X)$. This formula was obtained by using a probabilistic model which adjusted the local probabilities at each prime.
 
\begin{conj}[\textbf{Jones}]\label{jonesconj}
Let $E/\QQ$ be an elliptic curve without complex multiplication and let $L\geq 2$ be a positive integer. Then there is a non-negative real 
constant $C_{E,L}\geq 0$ such that, as $X\rightarrow\infty$, we have that
\begin{equation*}
\pi_{E,L}(X) \sim C_{E,L}\int_2^X\frac{1}{2\sqrt{t}(\log t)^L}dt.
\end{equation*}
\end{conj}

In Conjecture \ref{silverstangeamconj} we assume that there are infinitely many primes $p$ such that $\#E_p(\FF_p)$ is prime. Koblitz \cite{NK:1} gave the following conjecture for the number of primes $p\leq X$ such that $\#E_p(\FF_p)$ is prime, where the explicit constant in the asymptotic formula was refined by Zywina \cite{DZ:1}.

\begin{conj}[\textbf{Koblitz}]\label{koblitzconj}
Let $E/\QQ$ be an elliptic curve without complex multiplication. Then there exists a  
constant $C_{E}^{\rm twin}$ depending only on $E$ such that as $X\rightarrow\infty$
\begin{equation*}
\pi_{E}^{\rm twin}(X):=\#\{p\leq X: \#E_p(\FF_p) \:{\rm is\:prime}\} \sim C_{E}^{\rm twin}\frac{X}{(\log X)^2}.
\end{equation*}
\end{conj}

\begin{rmks} {\em (i) Jones \cite{NJ:1} showed that under the assumption of Conjecture \ref{koblitzconj} there are examples of elliptic curves such that $C_{E,L}=0$.

(ii) There are also other famous conjectures about the distributions of invariants associated with the reductions of elliptic curves over finite fields. These include the Sato-Tate conjecture for the distribution of the angles associated to the normalized traces $ \frac{a_p(E)}{2 \sqrt{p} }$ (we refer the reader to the survey paper \cite{KMRM:1} for an introduction) and the Lang-Trotter conjecture
\cite{SLHT:1} for the number of primes $p \leq X$ such that
$a_p(E) = t$ for a fixed integer $t$.

(iii) The Sato-Tate conjecture was recently proven for elliptic curves over totally real fields which have multiplicative 
reduction at some primes by Harris, Shepherd-Barron and Taylor \cite{HSBT:1}, but the other conjectures 
are completely open. For example, for the Lang-Trotter conjecture in the case $t\neq0$ we do not even know if there exist infinitely many primes $p$ such that
$a_p(E) = t$ for any elliptic curve over $\QQ$. The case $t=0$ corresponds to 
supersingular primes and was considered by Elkies \cite{NE:1}. He showed that every elliptic curve over $\QQ$ has infinitely many supersingular primes.
}
\end{rmks}

To gain insight into the above conjectures, it is natural to consider their averages over some family of elliptic 
curves. Let $a,b$ be integers and let $E_{a,b}$ be the elliptic curve given by the Weierstrass equation 
$$E_{a,b}: y^2=x^3+ax+b,$$ with the discriminant $\Delta(E_{a,b})\neq 0$. For $A,B>0$ we consider the two parameter family of elliptic curves
\begin{equation} \label{ecfamily}
\C:=\C(A,B)=\{E_{a,b}: |a|\leq A, |b|\leq B, \Delta(E_{a,b})\neq 0\}.
\end{equation}

In this paper we study the average for $\pi_{E,L}(X)$ over the family $\C(A,B)$ in $\eqref{ecfamily}$, that is, we consider the sum $\displaystyle{
 \frac{1}{|\C|}\sum_{E\in \C}\pi_{E,L}(X)}$. Our main result is the following theorem.

\begin{thrm} \label{aliquotupperbound}
 Let $\epsilon >0$, let $E/\QQ$ be an elliptic curve and let $\C$ be the family of elliptic curves in $\eqref{ecfamily}$ with
 $$A,B>X^{\epsilon} \quad {\rm and} \quad X^\frac{3L}{2}(\log X)^6< AB <e^{X^{\frac{1}{6}-\epsilon}}.$$ Then as $X\rightarrow \infty$ we have that $$
 \frac{1}{|\C|}\sum_{E\in \C}\pi_{E,L}(X)\ll_L \frac{\sqrt{X}}{(\log X)^L},$$ where the implied constant depends on $L$ only.
\end{thrm}

\begin{rmks}
{\em (i) Note that the additional condition $AB < e^{X^{\frac{1}{6}-\epsilon}}$ is not a limiting constraint since we are mainly 
interested in averages for small values of $A$ and $B$.\\
(ii) In $\eqref{alitriverror}$ we show that a trivial upper bound for the average is 
$$\frac{1}{|\C|}\sum_{E\in \C}\pi_{E,L}(X)\ll_L \sqrt{X}(\log\log X)^L$$ with 
$$A,B>X^L (\log X)^{L} (\log\log X)^L\quad {\rm and} \quad AB> X^{2L}(\log X)^{L} (\log\log X)^L.$$}
\end{rmks}

In Proposition \ref{hsievebound} we consider a sum of a product of class numbers over primes in a short interval. To obtain the conjectured upper bound for the average number of aliquot cycles over the family $\C$ we require the use of the fundamental lemma of sieve methods (see Lemma \ref{lemtwoone}) as well as a result of Granville and Soundararajan \cite{AGKS:1} (see Proposition \ref{lemtwothree}) to bound the error terms. This approach is also used in the work of Chandee, David, Koukoulopoulos and Smith \cite[Proposition 4.1]{CDKS:1}. However in their work, they are led to consider a sum of class numbers, whereas in our case we need to consider a sum of a product of class numbers.

To improve the bounds on $A$ and $B$, in Lemma \ref{errorimprove}, we consider the sum of aliquot cycles over representatives of isomorphism classes of elliptic curves. As in Banks and Shparlinski \cite{WBIS:1} and  Balog, Cojocaru, and David \cite{BCD:1}, we require the use of the large sieve inequality and a result of Friedlander and Iwaniec \cite{JFHI:2} (see Theorem \ref{fourthpower}). However, our calculations become much more technical since we must consider a product of $L$ characters.

\begin{rmks} {\em (i) Let $\epsilon>0$. The Lang-Trotter conjecture was shown to hold on average in the case $t=0$ for the family $\C(A,B)$ with $A,B>X^{\frac{1}{2}+\epsilon}$ and 
$AB>X^{\frac{3}{2}+\epsilon}$ by Fouvry and Murty \cite[Thoerem 6]{EFMM:1}. David and Pappalardi \cite{CDFP:1} then showed that the Lang-Trotter conjecture holds on average for any integer $t\neq 0$. The bounds on the size of $A$ and $B$ are an important feature of average results and several techniques for improving them have been developed. Baier \cite{SB:2} showed that the Lang-Trotter conjecture holds on average for any integer $t$ with $A,B > X^{\epsilon}$ and $AB > X^{3/2 + \epsilon}$. Banks and Shparlinski \cite{WBIS:1} used multiplicative 
character sums to show that the Sato-Tate Conjecture holds on average for the family ${\mathcal{C}}(A,B)$ with $A,B > X^{\epsilon}$ and
$AB > X^{1 + \epsilon}$. Finally, the Koblitz conjecture was shown to hold on average for the family $\C(A,B)$ with $A,B>X^{\epsilon}$ and 
$AB>X^{1+\epsilon}$ by Balog, Cojocaru, and David \cite{BCD:1}.}
\end{rmks}

Average results can give strong evidence for the distribution conjectures discussed above, because they also produce average conjectural constants in their respective asymptotic formulas. To derive a formula for the constant $C_{E,L}$ given in Conjecture 
\ref{jonesconj} we need to study ${\rm Prob}(\ell \nmid p+1-a_p(E))$ for primes $\ell$ and $p$.

For a non-zero integer $n$, we denote the $n$-torsion subgroup of $E$ by $E[n]$. Let 
$\QQ(E[n])$ be the field generated by adjoining to $\QQ$ the 
$x$ and $y$-coordinates of the $n$-torsion points of $E$. We have that $E[n]\cong \ZZ/n\ZZ \times \ZZ/n\ZZ$ for $n\geq 2$.
Since each element of the Galois group  ${\rm Gal}(\bar{\QQ}/\QQ)$ acts on $E[n]$ we have that ${\rm Gal}(\QQ(E[n])/\QQ)\subseteq {\rm GL}_2(\ZZ/n\ZZ)$ (see \cite[Chapter III.7]{JS:1}).

If $[{\rm GL}_2(\ZZ/n\ZZ):{\rm Gal}(\QQ(E[n])/\QQ)]\leq 2$ for each $n\geq 1$ 
(see \cite[pp. 309-311]{JPS:1} and \cite[p. 51]{SLHT:1}) then $E$ is called a \textit{Serre} curve. 
Jones \cite{NJ:1} has shown that for any Serre curve $E$, we have that $C_{E,L}>0$ and $C_{E,L}=C_L\cdot f_L(\Delta_{sf}(E))$, where $\Delta_{sf}(E)$ 
denotes the square-free part of the discriminant of any Weierstrass model of $E$ and $f_L$ is a positive function which approaches 1 
as $\Delta_{sf}(E)\rightarrow \infty$. In particular, for $L=2$, Jones \cite{NJ:1} gave the formula 
$$C_2=\frac{8}{3\pi^2}\prod_{\ell \:{\rm prime}} \frac{\ell^2(\ell^4-2\ell^3-2\ell^2+3\ell+3)}{((\ell^2-1)(\ell-1))^2}.$$ 

In a future work \cite{JP:1} we plan to verify the conjectural constant $C_2$ by obtaining an asymptotic result for the average of $\pi_{E,2}(X)$.

\subsection{Acknowledgment} This work constitutes a large portion of my PhD thesis. I thank my advisor, Chantal David for all her great 
advice and support while working on this problem. I would also like to thank Dimitris Koukoulopoulos and Amir Akbary for their helpful discussions related 
to this paper.

\section{Preliminaries}
For a basic introduction to the theory of elliptic curves we refer the reader to \cite{JS:1}.
Here, and in the rest of the paper, we let $\chi_d(n)$ denote the quadratic Dirichlet character defined by the Kronecker symbol namely, $$\chi_d(n):=\left(\frac{d}{n}\right).$$ We let $$L(s,\chi_d):=\sum_{n=1}^\infty \frac{\chi_d(n)}{n^s}= \prod_{\ell \:{\rm prime}}\left(1-\frac{\chi_d(\ell)}{\ell^s}\right)^{-1} \quad {\rm for\:Re}(s)>1,
$$ be the Dirichlet $L$-function associated to $\chi_d$. For $y>1$ we define the truncated quadratic Dirichlet $L$-function as $$L(1,\chi_d;y):=\prod_{\ell\leq y}\left(1-\frac{\chi_d(\ell)}{\ell}\right)^{-1}.$$ 

The following proposition is a consequence of a result of Granville and Soundararajan \cite{AGKS:1} essentially due to Elliot \cite{PE:1}. It allows us to bound the error terms in our calculations in Proposition \ref{hsievebound}.

\begin{prop}[\textbf{Granville-Soundararajan}] \label{lemtwothree}
Let $\alpha\geq 1$ and $Q\geq 3$. There is a set ${\mathcal{E}}_\alpha(Q)\subset[1,Q]$ of at most 
$Q^{\frac{2}{\alpha}}$ integers such that if $\chi$ is a quadratic Dirichlet character of conductor $q\leq Q$ not in ${\mathcal{E}}_\alpha(Q)$, then 
$$L(1,\chi)=L(1,\chi;(\log Q)^{8\alpha^2})\left(1+O_\alpha\left(\frac{1}{(\log Q)^\alpha}\right)\right).$$
\end{prop}

\begin{proof}
The result is stated in terms of primitive characters in \cite[Proposition 2.2]{AGKS:1}. The proof of the proposition in its present form is given in \cite[Lemma 2.2]{CDKS:1}. 
\end{proof}

We now state the analytic class number formula for quadratic Dirichlet $L$-functions, (see Davenport \cite[Chapter 6]{HD:1}).

\begin{thrm} 
Let $D= df^2$ be a negative number such that $d$ is a negative fundamental discriminant and let $\chi_D$ be the Kronecker symbol. Then 
$$\frac{h(d)}{w(d)}= \frac{\sqrt{-D}}{2 \pi} L(1,\chi_D)$$ where $h(d)$ denotes the usual class number of the imaginary quadratic order of discriminant $d$ and $w(d)$ is the number of roots of unity in $\QQ(\sqrt{d})$.
\end{thrm}

We recall the following formulation of the definition of the Hurwitz-Kronecker class number, (see Lenstra \cite{HL:1}). Let $D$ be a negative (not necessarily fundamental) discriminant then the \textit{Hurwitz-Kronecker class number} of discriminant $D$ is defined by
\begin{equation*}
H(D)=\sum_{\substack{f^2 | D \\ \frac{D}{f^2}\equiv 0,1 \imod{4}}}\frac{h\left(\frac{D}{f^2}\right)}{w\left(\frac{D}{f^2}\right)}. 
\end{equation*}
This leads to the following useful result of Deuring \cite{MD:1}.

\begin{thrm}[\textbf{Deuring}] \label{deuring}
Let $p>3$ be a prime and let $t$ be an integer such that $t^2-4p<0$. Then
$$\sum_{\substack{\bar{E}/\FF_{p} \\ a_p(\bar{E})=t}} \frac{1}{\#{\rm Aut}(\bar{E})}=H(t^2-4p),$$
where $\bar{E}$ denotes a representative of an isomorphism class of $E/\FF_p$. 
\end{thrm}
As in the proof of Balog, Cojocaru, and David \cite[Lemma 6]{BCD:1} we require the following two theorems in the proof of Lemma \ref{errorimprove}. We first state the large sieve inequality for Dirichlet characters, for a proof, we refer the reader to Davenport \cite[Chapter 27]{HD:1}.

\begin{thrm} \label{largesieve}
Let $M, N, Q$ be positive integers and let $\left\{a_n\right\}_n$ be a sequence of complex numbers. For a fixed $q \leq Q$, we let 
$\chi$ be a Dirichlet character modulo $q$. Then
$$\sum_{q\leq Q}\frac{q}{\phi(q)}\sum_{\substack{\chi \imod{q} \\ \chi\:{\rm primitive}}}\left|\sum_{M< n\leq M+N} a_n\chi(n)\right|^2
 \leq (N+3Q^2)\sum_{M< n\leq M+N}|a_n|^2.$$
\end{thrm}
The second theorem is a result of Friedlander and Iwaniec \cite{JFHI:2} that bounds the fourth power moment of Dirichlet characters.

\begin{thrm}[\textbf{Friedlander-Iwaniec}] \label{fourthpower}
Let $q$ and $N$ be positive integers. Let $\chi$ denote a Dirichlet character modulo $q$, with $\chi_0$ denoting the principal character. 
Then $$ \sum_{\chi\neq \chi_0}\left|\sum_{n\leq N}\chi(n)\right|^4 \ll N^2q\log^6 q.$$
\end{thrm}

Finally, we end this section with a result known as the fundamental lemma of sieve methods. It is stated in various forms in the literature (see Halberstam and Richert \cite[p. 82]{HHHR:1} and Iwaniec and Kowalski \cite[Lemma 6.3]{HIEK:1}). The version we will use is a direct consequence of \cite[Lemma 5]{JFHI:1}. Here and throughout the rest of the paper we let $P^+(n)$ denote the largest prime dividing $n$ and let $P^-(n)$ denote the smallest prime dividing $n$. We denote 
by $(f*g)(n)$ the convolution $$(f*g)(n):=\sum_{d\mid n} f(d)g\left(\frac{n}{d}\right).$$

\begin{lem}\label{lemtwoone} 
Let $y\geq 2, D=y^{u}$ with $u\geq 2.$ There exists two arithmetic functions $\lambda^{\pm}\: :\NN\rightarrow [-1,1],$ 
supported in the set $\{d\in\NN: P^+(d) \leq y, d\leq D\}$, for which
$$\begin{cases}
(\lambda^-*1)(n)=(\lambda^+*1)(n)=1&\:{\rm if}\:P^-(n)>y, \\
(\lambda^-*1)(n)\leq 0 \leq (\lambda^+*1)(n) &\:{\rm otherwise}. \end{cases}$$
Moreover, if $g: \NN\rightarrow\RR$ is a multiplicative function with $0\leq g(p)\leq \min\{2,p-1\}$ for all primes $p\leq y$ then
$$\sum_{d}\frac{\lambda^\pm(d)g(d)}{d}=\prod_{p\leq y}\left(1-\frac{g(p)}{p}\right)(1+O(e^{-u})).$$
\end{lem}

\section{Reduction to an average of class numbers}
In this section we prove the main result, Theorem \ref{aliquotupperbound}. We begin this section by fixing notational conventions that we use for the remainder of the paper.

Let $P:=(p_1,\ldots,p_L)$ be a vector of $L$ distinct primes and denote the smallest prime in the vector as $p:=p_{L+1}:=p_1$.  
For a fixed elliptic curve $E_{a,b}$, we define the following indicator function which determines if $P$ is a normalized aliquot cycle of 
length $L$, 
\begin{equation*}
\wpab:=\begin{cases} 1 & {\rm if}\: \#E_{p_i,a,b}(\FF_{p_i})=p_{i+1}\:{\rm for}\:1\leq i\leq L, \\ 0 & {\rm otherwise}. \end{cases}
\end{equation*}
Let $S:=(s_1,\ldots,s_L)$ and $T:=(t_1,\ldots,t_L)$ be vectors such that $s_i,t_i\in\FF_{p_i}$ for $1\leq i\leq L$. This leads to the similar function,
\begin{equation} \label{defnofwpst}
w(P,S,T):=\begin{cases} 1 & {\rm if}\: \#E_{p_i,s_i,t_i}(\FF_{p_i})=p_{i+1}\:{\rm for}\:1\leq i\leq L, \\ 0 & {\rm otherwise}. \end{cases}
\end{equation}

We also define the following products 
\begin{equation*}
\FF(P):=\FF_{p_1}\times\cdots\times\FF_{p_L} \quad {\rm and}\quad
\FF(P)^*:=\FF_{p_1}^*\times\cdots\times\FF_{p_L}^*.
\end{equation*}
Thus,
\begin{equation*}
\sum_{\FFp}1=\sum_{\substack{1\leq s_1\leq p_1 \\ 1\leq t_1\leq p_1}}\cdots\sum_{\substack{1\leq s_L\leq p_L \\ 
1\leq t_L\leq p_L}}1 \quad{\rm and} \quad \sum_{\FFps}1=\sum_{\substack{1\leq s_1< p_1 \\ 
1\leq t_1< p_1}}\cdots\sum_{\substack{1\leq s_L< p_L \\ 1\leq t_L< p_L}}1.
\end{equation*}
For positive integers $m$ and $n$ we define the symmetric function that arises from the 
application of Theorem \ref{deuring} $$D(m,n):=(m+1-n)^2-4m=(n+1-m)^2-4n=D(n,m).$$ 

Finally, we recall the definitions of $\eqref{pplusminus}$ and $\eqref{ecfamily}$. We denote the sum 
over $P$ as 
$$\sumofallps 1:= \sum_{p_1\leq X}\sum_{p_1^-<p_2<p_1^+}\cdots \sum_{p_{L-1}^-< p_{L} < p_{L-1}^+}1,$$
and we have that $|\C|=4AB + O(A+B+1)$.

We begin by considering the trivial upper bound for the average number of aliquot cycles. We have that
\begin{align}
&\frac{1}{|\C|} \sum_{E\in \C} \pi_{E,L}(X)\nonumber\\
=&\frac{1}{|\C|}\sum_{E_{a,b}\in \C}\sumofallps \wpab=\frac{1}{|\C|}\sumofallps\sum_{E_{a,b}\in \C}\wpab\label{alibegin}\\
=&\frac{1}{|\C|}\sumofallps\sum_{\FFp}w(P,S,T)\sum_{\substack{|a|\leq A, |b|\leq B \\ a\equiv s_i 
\imod{p_i} \\ b\equiv t_i \imod{p_i} \\ 1\leq i \leq L}}1\nonumber\\
=&\frac{4AB}{|\C|}\sumofallps\sum_{\FFp}\frac{w(P,S,T)}{p_1^2\cdots p_L^2}+O\Bigg(\frac{(B+A)}{|\C|}\sumofallps\sum_{\FFp}\frac{w(P,S,T)}{p_1\cdots p_L}\nonumber\\
+&\frac{1}{|\C|}\sumofallps\sum_{\FFp}w(P,S,T)\Bigg),\label{alitrivone}
\end{align}
where
\begin{equation}
\sum_{\FFp}w(P,S,T)=\sum_{\substack{1\leq s_1,t_1\leq p_1 \\ \#E_{p_1,s_1,t_1}(\FF_{p_1})=p_2}}\cdots\sum_{\substack{1\leq s_L,t_L\leq p_L \\ \#E_{p_L,s_L,t_L}(\FF_{p_L})=p_1}}1. \label{wexpresstriv}
\end{equation}
For $1\leq i\leq L$ the sums in $\eqref{wexpresstriv}$ over $s_i$ and $t_i$ can be changed to a sum over isomorphism classes which we 
denote by $\bar{E}_{p_i,s_i,t_i}$. Then we have that
\begin{align}
\sum_{\substack{1\leq s_i,t_i\leq p_i \\ \#E_{p_i,s_i,t_i}(\FF_{p_i})={p_{i+1}}}}1&=\sum_{\substack{\bar{E}_{p_i,s_i,t_i}/ \FF_{p_i} \\ p_i+1-a_{p_i}(\bar{E}_{p_i,s_i,t_i})=p_{i+1}}}\frac{p_i-1}{\#\aut(\bar{E}_{p_i,s_i,t_i}(\FF_{p_i}))}\nonumber \\
&=(p_i-1)H((p_i+1-p_{i+1})^2-4p_i)=(p_i-1)H(D(p_i,p_{i+1})), \label{alitrivtwo}
\end{align}
by Theorem \ref{deuring}. From the convexity bound for a Dirichlet character $\chi$ of modulus $d$, we have that $L(1,\chi_d)\ll\log |d|$. Therefore, by the 
analytic class number formula for $1\leq i\leq L$, we deduce that
\begin{align}
H(D(p_i,p_{i+1}))&=\sum_{\substack{f^2| D(p_i,p_{i+1}) \\ \frac{ D(p_i,p_{i+1})}{f^2}\equiv 0,1 \imod{4}}} \frac{\sqrt{| D(p_i,p_{i+1})|}}{2\pi f}L\left(1,\left(\frac{ D(p_i,p_{i+1})/f^2}{\cdot}\right)\right)\nonumber\\
&\ll\sqrt{ |D(p_i,p_{i+1})|}(\log p_i)\sum_{f\mid  D(p_i,p_{i+1})}\frac{1}{f}\ll \sqrt{p_i}(\log p_i )(\log\log |D(p_i,p_{i+1})|)\nonumber \\
&\ll \sqrt{p}(\log p )(\log \log p), \label{hdellbound}
\end{align}
since $p_i=p+O(\sqrt{p})$.

Thus, from $\eqref{alitrivtwo}$ and $\eqref{hdellbound}$ we have that the main term in $\eqref{alitrivone}$ is bounded by
\begin{align}
&\frac{AB}{|\C|}\sum_{p\leq X}\frac{1}{p^L}\sumofallh\nonumber\\
\ll_L&\left(1+O\left(\frac{1}{A}+\frac{1}{B}+\frac{1}{AB}\right)\right)\sum_{p\leq X}\frac{1}{p^L}\frac{p^{\frac{L-1}{2}}}{(\log p)^{L-1}}p^\frac{L}{2}(\log p)^{L}(\log \log p)^{L}\nonumber\\
\ll_L&\sum_{p\leq X}\frac{\log p(\log \log p)^L}{\sqrt{p}}\ll_L \sqrt{X}(\log\log X)^L.\label{alitrivmain}
\end{align}
Similarly, the error term in $\eqref{alitrivone}$ is bounded by
\begin{align}
&\left(\frac{1}{A}+\frac{1}{B}\right)X^{L+\frac{1}{2}}(\log\log X)^L+\frac{X^{2L+\frac{1}{2}}(\log\log X)^L}{AB}.\label{alitriverror}
\end{align}
Hence, from $\eqref{alitriverror}$ to obtain the correct upper bound for the average we need $$A,B>X^L (\log X)^{L} (\log\log X)^L\quad {\rm and} \quad AB> X^{2L}(\log X)^{L} (\log\log X)^L,$$
whereas $\pi_{E,L}(X)$ only considers primes of size at most $X$. 
Also, we see that using the bound from $\eqref{hdellbound}$ for $H(D(p_i,p_{i+1}))$ in $\eqref{alitrivone}$ 
does not give the correct order of magnitude for the main term in $\eqref{alitrivmain}$. 
Therefore, to obtain the conjectured upper bound for Theorem \ref{aliquotupperbound} we develop techniques not present in the estimations above. This is the approach of the following theorem.

\begin{thrm} \label{aliquotaverage}
Let $\epsilon >0$, let $E/\QQ$ be an elliptic curve and let $\C$ be the family of elliptic curves in $\eqref{ecfamily}$ with
 $$A,B>X^{\epsilon} \quad {\rm and} \quad X^\frac{3L}{2}(\log X)^6< AB<e^{X^{\frac{1}{6}-\epsilon}}.$$ Then as $X\rightarrow \infty$ we have that 
\begin{equation}
\frac{1}{|\C|}\sum_{E\in \C}\pi_{E,L}(X)=\Bigg(\sumofallps \prod_{j=1}^L\frac{H(D(p_j,p_{j+1}))}{p_j}\Bigg)\left(1+
O\left(\frac{1}{X^\epsilon}\right)\right).\label{avgthrmsum}
\end{equation}
\end{thrm}

We have that the sum on the RHS of $\eqref{avgthrmsum}$ is
\begin{align*}
&\sum_{p\leq X}\frac{1}{p^L}\prod_{i=1}^{L-2}\Bigg(\sum_{p_i^-<p_{i+1}<p_i^+}H(D(p_i,p_{i+1}))\Bigg)
\sum_{p_{L-1}^-<p_{L}< p_{L-1}^+}H(D(p_{L-1},p_L))H(D(p_L,p))\\
\times&\left(1+O_L\left(\frac{1}{\sqrt{p}}\right)\right),
\end{align*} 
since $p_i= p+O(\sqrt{p})$ for $1< i \leq L$. We use the following technical propositions to bound the inner sums above. 

\begin{prop} \label{hsievebound}
Fix primes $p, r>3$ not necessarily distinct with $r=p+O(\sqrt{p})$ and let $q$ be a prime in the range $p^-< q < p^+$ with $q\neq p$ or $r$. Then we have that
\begin{equation*}
\sum_{p^-< q < p^+} H(D(p,q))H(D(r,q)) \ll \frac{p^\frac{3}{2}}{\log p}.
\end{equation*}
\end{prop}

\begin{prop} \label{hsinglebound}
Let $p$ and $q$ be distinct primes such that $p^-< q < p^+$. Then we have that
\begin{equation*}
\sum_{p^-< q < p^+} H(D(p,q)) \ll \frac{p}{\log p}.
\end{equation*}
\end{prop}

We delay the proofs of Proposition \ref{hsievebound} and Proposition \ref{hsinglebound} until the following section. 
We now have that Theorem \ref{aliquotupperbound} is an immediate consequence of Theorem \ref{aliquotaverage}.

\begin{proof} (Proof of Theorem \ref{aliquotupperbound})
From Proposition \ref{hsievebound} and Proposition \ref{hsinglebound} we have by partial summation that the main 
term in $\eqref{avgthrmsum}$ is
\begin{align*}
\sumofallps \prod_{j=1}^L\frac{H(D(p_j,p_{j+1}))}{p_j}=&\sum_{p\leq X}\frac{1}{p^L}\prod_{i=1}^{L-2}
\Bigg(\sum_{p_i^-<p_{i+1}<p_i^+}H(D(p_i,p_{i+1}))\Bigg)\\
\times&\sum_{p_{L-1}^-<p_{L}< p_{L-1}^+}H(D(p_{L-1},p_L))H(D(p_L,p))\left(1+O_L\left(\frac{1}{\sqrt{p}}\right)\right)\\
\ll_L&\sum_{p\leq X}\frac{1}{p^L}\frac{p^{L-2}}{(\log p)^{L-2}}\frac{p^\frac{3}{2}}{\log p}=\sum_{p\leq X}
\frac{1}{\sqrt{p}(\log p)^{L-1}} \ll_L \frac{\sqrt{X}}{(\log X)^L}.
\end{align*}
\end{proof}

\begin{proof} (Proof of Theorem \ref{aliquotaverage})
 We begin the proof by recalling $\eqref{alibegin}$,  $$\frac{1}{|\C|} \sum_{E\in \C} \pi_{E,L}(X)=\frac{1}{|\C|}
\sumofallps\sum_{E_{a,b}\in \C}\wpab.$$ To obtain an improvement on this sum, instead of summing over elliptic curves, we will sum over representatives of isomorphism classes. Let $E_{s,t}$ be an elliptic curve defined over $\FF_{p}$. We count the curves 
$E_{a,b} \in \C$ whose reductions modulo $p$ are isomorphic to $E_{s,t}$ over $\FF_{p}$. Recall that two elliptic curves 
$E_{s,t}$ and $E_{s',t'}$ are isomorphic over $\FF_{p}$ if and only if there exists a $u\in \FF_{p}^*$ such that $s'=su^4$ and 
$t'=tu^6$. Thus, we have that the number of elliptic curves over $\FF_{p}$ isomorphic to $E_{s,t}$ is 
$$\frac{\#\FF_{p}^*}{\#\aut(E_{s,t})}= \frac{p-1}{\#\aut(E_{s,t})}.$$ More precisely, if we are counting $|a|\leq A, |b|\leq B$ such that if there 
exists $u_i\in \FF_{p_i}^*$ such that $a\equiv s_iu_i^4 \imod{p_i}$ and $b\equiv t_iu_i^6 \imod{p_i}$ then for each fixed elliptic curve 
$E_{s_i,t_i}$ we will be over counting by the number of elliptic curves over $\FF_{p_i}$ isomorphic to $E_{s_i,t_i}$. By correcting for 
this over count we have that the sum over elliptic curves in $\eqref{alibegin}$ 
becomes
$$\sum_{E_{a,b}\in \C}\wpab= \sum_{\FFp}w(P,S,T)\prod_{j=1}^L\frac{\#\aut(E_{p_j,s_j,t_j})}{(p_j-1)}
\sum_{\substack{|a|\leq A, |b|\leq B \\ \exists (u_1,\ldots,u_L) \in \FF(P)^*\\ a \equiv s_iu_i^4 \imod{p_i},  
b \equiv t_iu_i^6 \imod{p_i} \\ 1\leq i \leq L}}1.$$ 
Hence, $\eqref{alibegin}$ becomes
\begin{equation}
\frac{1}{|\C|} \sum_{E\in \C} \pi_{E,L}(X)=\frac{1}{|\C|}\sumofallps\sum_{\FFp}w(P,S,T)R(P,S,T)
\prod_{j=1}^L\frac{\#\aut(E_{p_j,s_j,t_j})}{(p_j-1)}, \label{aliquotaut}
\end{equation}
where $R(P,S,T)$ is the number of integers $|a|\leq A,|b|\leq B$ such that there exists a vector
$(u_1,\ldots,u_L) \in \FF(P)^*$ satisfying
\begin{equation} \label{rpst}
 a\equiv s_iu_i^4\imod{p_i}, \quad b\equiv t_iu_i^6\imod{p_i} 
\quad{\rm for}\:1\leq i\leq L.
\end{equation}

For an elliptic curve $E_{s,t}/\FF_p$, we have that the order of the automorphism group of $E_{s,t}$ is given by
\begin{equation*}
\#\aut(E_{s,t})=\begin{cases} 6 & {\rm if}\: s=0\:{\rm and}\: p\equiv 1 \imod{3}, \\
4 & {\rm if}\: t=0\:{\rm and}\: p\equiv 1 \imod{4}, \\
2 & {\rm otherwise}.
\end{cases}
\end{equation*}
Thus, we split up the sum in $\eqref{aliquotaut}$ into two cases, $s_it_i \neq 0$ and $s_it_i = 0$ 
to write $\eqref{aliquotaut}$ as
\begin{align}
\frac{1}{|\C|} \sum_{E\in \C} \pi_{E,L}(X)=&\frac{2^L}{|\C|}\sumofallps\sum_{\FFps}
\frac{w(P,S,T)R(P,S,T)}{(p_1-1)\cdots(p_L-1)}
\nonumber \\
+&\frac{1}{|\C|}\sumofallps\sum_{\substack{\FFp \\ s_it_i= 0 \\{\rm for\:some}\: 1\leq i \leq L}}w(P,S,T)R(P,S,T)
\prod_{j=1}^L\frac{\#\aut(E_{p_j,s_j,t_j})}{(p_j-1)}.\label{alione}
\end{align}

We can express the first sum in $\eqref{alione}$ as
\begin{align}
&\frac{4AB}{|\C|}\sumofallps\prod_{j=1}^L\frac{1}{p_j(p_j-1)}\sum_{\FFps}w(P,S,T) \nonumber\\
+&\frac{2^L}{|\C|}\sumofallps\prod_{j=1}^L\frac{1}{(p_j-1)}\sum_{\FFps}w(P,S,T)\left(R(P,S,T)-
\frac{4AB}{2^L p_1\cdots p_L}\right).\label{alitwo}
\end{align}
The first term in $\eqref{alitwo}$ contributes to the main term and we use the following technical lemma, where we delay its proof to Section 5, to bound the second
term in $\eqref{alitwo}$.

\begin{lem} \label{errorimprove}
Let $L\geq 2$ be an integer, let $E/\QQ$ be an elliptic curve and let $A,B>0$. Then for any positive integer $k$, as $X\rightarrow \infty$ we have that
\begin{align}
&\sumofpfracsing\sum_{\FFps}w(P,S,T) \left(R(P,S,T)-\frac{AB}{2^{L-2}p_1\cdots p_L}\right)\nonumber \\
\ll_{k,L} &ABX^{\frac{1}{2}-\frac{L+1}{4k}}(\log X)^\frac{L}{2k}(\log\log X)^L\left((\log A)^{\frac{k^2-1}{2k}}+(\log B)^{\frac{k^2-1}{2k}}\right)\nonumber\\
+&(A\sqrt{B}+B\sqrt{A})X^{\frac{1}{2}+\frac{3L-1}{4k}}(\log X)^{\frac{k^2+L-1}{2k}}(\log\log X)^L+
\sqrt{AB}X^{\frac{3L+2}{4}}(\log X)^{3-L},\label{lemerrorbound}
\end{align}
where $w(P,S,T)$ is given in $\eqref{defnofwpst}$ and $R(P,S,T)$ is given in $\eqref{rpst}$.
\end{lem}
Thus, from Lemma \ref{errorimprove} we have that for any positive integer $k$, the second sum in $\eqref{alitwo}$ becomes
\begin{align}
&\frac{2^L}{|\C|}\sumofpminusone\sum_{\FFps}w(P,S,T)\left(R(P,S,T)-\frac{4AB}{2^L p_1\cdots p_L}\right)\nonumber \\
\ll_{L,k}&X^{\frac{1}{2}-\frac{L+1}{4k}}(\log X)^\frac{L}{2k}(\log\log X)^L\left((\log A)^{\frac{k^2-1}{2k}}+(\log B)^{\frac{k^2-1}{2k}}\right) +\frac{1}{\sqrt{AB}}
X^{\frac{3L+2}{4}}(\log X)^{3-L}\nonumber \\
+& \left(\frac{1}{\sqrt{B}}+\frac{1}{\sqrt{A}}\right)X^{\frac{1}{2}+\frac{3L-1}{4k}}(\log X)^{\frac{k^2+L-1}{2k}}(\log\log X)^L. 
\label{aliquotlembound}
\end{align}

We now consider the inner sum in the first sum in $\eqref{alitwo}$,
\begin{equation}
\sum_{\FFps}w(P,S,T)=\sum_{\substack{1\leq s_1,t_1<p_1 \\ \#E_{p_1,s_1,t_1}(\FF_{p_1})=p_2}}\cdots\sum_{\substack{1\leq s_L,t_L<p_L \\ 
\#E_{p_L,s_L,t_L}(\FF_{p_L})=p_1}}1. \label{wexpress}
\end{equation}
Similarly to the calculation of $\eqref{alitrivtwo}$ we have by Theorem \ref{deuring} that 
\begin{align}
\sum_{\substack{1\leq s_i,t_i<p_i \\ \#E_{p_i,s_i,t_i}(\FF_{p_i})={p_{i+1}}}}1
&=\sum_{\substack{\bar{E}_{p_i,s_i,t_i}/\FF_{p_i} \\ p_i+1-a_{p_i}(\bar{E}_{p_i,s_i,t_i})=p_{i+1}}}\frac{p_i-1}
{\#\aut(\bar{E}_{p_i,s_i,t_i}(\FF_{p_i}))}+O(p_i) \nonumber \\
&=(p_i-1)H((p_i+1-p_{i+1})^2-4p_i)+O(p_i). \label{alithree}
\end{align}
Thus, from $\eqref{hdellbound}$ and $\eqref{alithree}$ we have that $\eqref{wexpress}$ becomes
\begin{align}
\sum_{\FFps}w(P,S,T)=&\prod_{i=1}^{L}\left((p_i-1)H(D(p_i,p_{i+1}))+O(p_i)\right) \nonumber\\
=&\prod_{i=1}^{L}(p_i-1)H(D(p_i,p_{i+1}))+O_L\left(p^{\frac{3L-1}{2}}(\log p)^{L-1}(\log\log p)^{L-1}\right).\label{dubsum}
\end{align}
Combining $\eqref{dubsum}$ with the first term in $\eqref{alitwo}$ gives
\begin{align}
&\frac{4AB}{|\C|}\sumofallps\prod_{j=1}^L\frac{1}{p_j(p_j-1)}\sum_{\FFps}w(P,S,T)\nonumber\\
=&\frac{4AB}{|\C|}\sumofallps\Bigg(\prod_{j=1}^{L}\frac{H(D(p_j,p_{j+1}))}{p_j}+O_L\left(\frac{1}{p^{2L}} 
\cdot p^{\frac{3L-1}{2}}(\log p)^{L-1}(\log\log p)^{L-1}\right)\Bigg) \nonumber\\
=&\Bigg(\sumofallps\prod_{j=1}^L\frac{H(D(p_j,p_{j+1}))}{p_j}+O_L\left((\log\log X)^{L}\right)\Bigg)\left(1+O_L\left(\frac{1}{A}+\frac{1}{B}+\frac{1}{AB}\right)\right).\label{alifive}
\end{align}

We see that the first term in $\eqref{alifive}$ gives the main term in $\eqref{avgthrmsum}$ and 
by Proposition \ref{hsievebound} and Proposition \ref{hsinglebound} we have that the error term in $\eqref{alifive}$ is 
bounded by
\begin{align*}
 &\left(\sum_{p\leq X}\frac{1}{p^L}\frac{p^{L-2}}{(\log p)^{L-2}}\frac{p^{\frac{3}{2}}}{\log p}\right)\left(\frac{1}{A}+
\frac{1}{B}+\frac{1}{AB}\right)+(\log\log X)^{L}) \\
\ll_L&\left(\frac{1}{A}+\frac{1}{B}+\frac{1}{AB}\right)\sum_{p\leq X}\frac{1}{\sqrt{p}(\log p)^{L-1}}\ll_L\left(\frac{1}{A}+\frac{1}{B}+\frac{1}{AB}\right)\frac{\sqrt{X}}{\log X},
\end{align*}
which is smaller than the second and third terms in the error terms in $\eqref{aliquotlembound}$.

Thus, it remains to consider the second term in $\eqref{alione}$. Similarly to the treatment of the average of the Lang-Trotter 
Conjecture by Baier \cite[Theorem 2.1]{SB:2} we have that
\begin{align}
&\frac{1}{|\C|}\sumofallps\sum_{\substack{\FFp \\ s_it_i= 0  \\{\rm for\:some}\: 1\leq i \leq L}}w(P,S,T)
\prod_{j=1}^L\frac{\#\aut(E_{p_j,s_j,t_j})}{(p_j-1)}\sum_{\substack{|a|\leq A, |b|\leq B \\ \FFpu \\ a\equiv s_iu_i^4 \imod{p_i} \\ 
b\equiv t_iu_i^6 \imod{p_i}\\ 1\leq i \leq L}}1\nonumber \\
\ll_L&\frac{1}{|\C|}\sumofallps\sum_{\substack{|a|\leq A, |b|\leq B \\ ab\equiv 0 \imod{p_1} \:{\rm or } \\ ab\equiv 0 \imod{p_i} 
\\ {\rm for} \:2\leq i \leq L}}\wpab. \label{baierone}
\end{align}
If $ab\equiv 0 \imod{p_j}$ then fixing $p_j$ completely determines the other $p_i$ for $1\leq i\neq j \leq L$ from $\wpab$. Hence, without loss of generality we can assume that  $ab\equiv 0 \imod{p_1}$ and we have that $\eqref{baierone}$ is bounded by
\begin{equation}
 \frac{1}{|\C|}\sum_{\substack{|a|\leq A \\  |b|\leq B}}\sum_{\substack{p\leq X \\ p \mid ab}}\wpab
\ll_L\frac{1}{|\C|}\sum_{\substack{|a|\leq A, |b|\leq B }}\tau(ab)\ll_L\frac{1}{|\C|}\sum_{n\leq AB}\tau^2(n)\ll_L (\log AB)^3. \label{baierbound}
\end{equation}
From $\eqref{aliquotlembound}, \eqref{alifive}$ and $\eqref{baierbound}$ we have that
\begin{align}
\frac{1}{|\C|}\sum_{E\in \C}&\pi_{E,L}(X)=\sum_{p\leq x}\frac{1}{p^{L}}\sumofallh+O_{L,k}\Bigg((\log AB)^3\nonumber\\
+&X^{\frac{1}{2}-\frac{L+1}{4k}}(\log X)^\frac{L}{2k}(\log\log X)^L\left((\log A)^{\frac{k^2-1}{2k}}+(\log B)^{\frac{k^2-1}{2k}} \right)+\frac{1}{\sqrt{AB}}X^{\frac{3L+2}{4}}(\log X)^{3-L} \nonumber \\
+& \left(\frac{1}{\sqrt{B}}+\frac{1}{\sqrt{A}}\right)X^{\frac{1}{2}+\frac{3L-1}{4k}}(\log X)^{\frac{k^2+L-1}{2k}}(\log\log X)^L\Bigg). 
\label{bigaliquott}
\end{align}

Now the first term in the error term of $\eqref{bigaliquott}$ is smaller than the main term if 
$$AB<e^\frac{X^{1/6}}{(\log X)^{L/3}}.$$
The second term in the error term of $\eqref{bigaliquott}$ is smaller than the main term for any $k\geq 1$.
The third term in the error term of $\eqref{bigaliquott}$ is smaller than the main term if $$AB> X^\frac{3L}{2}(\log X)^{6}.$$
The fourth term in the error term of $\eqref{bigaliquott}$ is smaller than the main term if 
$$A,B > X^{\frac{3L-1}{2k}}(\log X)^{\frac{k^2+L-1}{k}+2L}(\log\log X)^{2L}.$$
For every $\epsilon>0$ we can find a positive integer $k$ such that $$\epsilon > \frac{3L-1}{2k},$$ and therefore the fourth term 
in the error term of $\eqref{bigaliquott}$ is smaller than the main term if $A,B>X^\epsilon$, which gives the result.
\end{proof}

\section{Upper bounds on sums of class numbers}
\begin{proof} (Proof of Proposition \ref{hsievebound}) 
We begin by using the analytic class number formula to relate the class number $H(D)$ to a quadratic Dirichlet 
$L$-function evaluated at one. We have that
\begin{align*}
\sum_{p^-< q < p^+} H(D(p,q))H(D(r,q))=&\sum_{p^-< q < p^+}\sum_{\substack{f_1^2|D(p,q) \\ 
(f_1,2)=1}} \frac{\sqrt{|D(p,q)|}}{2\pi f_1}L\left(1,\left(\frac{D(p,q)/f_1^2}{\cdot}\right)\right)\\
\times&\sum_{\substack{f_2^2|D(r,q) \\ (f_2,2)=1}} \frac{\sqrt{|D(r,q)|}}{2\pi f_2}
L\left(1,\left(\frac{D(r,q)/f_2^2}{\cdot}\right)\right),
\end{align*}
since $\frac{D(p,q)}{f^2} \not \equiv 0 \imod{4}$ for $p,q>3$ and $\frac{D(p,q)}{f^2}\equiv 1 \imod{4}$ if and only if $f$ is odd. We also have that $q,r=p+O(\sqrt{p})$ and hence, $D(p,q), D(r,q) \ll p$. With the goal of obtaining an upper bound for the LHS of the above identity we define the sum
\begin{equation} \label{mqprop}
S_1:=\sum_{p^-< q < p^+}\sum_{\substack{f_1^2|D(p,q) \\ f_2^2|D(r,q) \\ (f_1f_2,2)=1}}
\frac{L\left(1,\left(\frac{D(p,q)/f_1^2}{\cdot}\right)\right)L\left(1,\left(\frac{D(r,q)/f_2^2}{\cdot}\right)\right)}{f_1f_2}.
\end{equation}

We have that
\begin{align}
L\left(1,\left(\frac{D(p,q)/f_1^2}{\cdot}\right)\right) 
&= \prod_\ell \left(1- \left(\frac{D(p,q)/f_1^2}{\ell}\right)\frac{1}{\ell}\right)^{-1}
\leq \frac{2f_1}{\varphi(f_1)}\prod_{\ell \nmid 2f_1}\left(1-\frac{\left(\frac{(2f_1)^2D(p,q)}{\ell}\right)}{\ell}\right)^{-1}\nonumber\\
&\ll\frac{f_1}{\varphi(f_1)}L\left(1,\left(\frac{(2f_1)^2D(p,q)}{\cdot}\right)\right),\label{lwithfs}
\end{align}
and similarly, $$L\left(1,\left(\frac{D(r,q)/f_2^2}{\cdot}\right)\right)\ll \frac{f_2}{\varphi(f_2)}
L\left(1,\left(\frac{(2f_2)^2D(r,q)}{\cdot}\right)\right).$$

To ease notation for the remainder of this section we denote $$\chi_1:=\left(\frac{(2f_1)^2D(p,q)}{\cdot}\right)
 \quad {\rm and} \quad \chi_2:=\left(\frac{(2f_2)^2D(r,q)}{\cdot}\right).$$ Now we have that
\begin{equation}
S_1\ll \sum_{p^-< q < p^+} \sum_{\substack{f_1^2|D(p,q) \\ f_2^2|D(r,q) \\ (f_1f_2,2)=1}}\frac{L\left(1,\chi_1\right)
L\left(1,\chi_2\right)}{\varphi(f_1)\varphi(f_2)}\ll \sum_{p^-< q < p^+} \sum_{\substack{f_1|D(p,q) \\ f_2|D(r,q) \\ 
(f_1f_2,2)=1}}\frac{L\left(1,\chi_1\right)L\left(1,\chi_2\right)}{\varphi(f_1)\varphi(f_2)},\label{mqbound}
\end{equation}
since the sum on the RHS in $\eqref{mqbound}$ is larger than the sum in $\eqref{mqprop}$. Then $$\sum_{p^-< q < p^+} 
H(D(p,q))H(D(r,q)) \ll pS_1.$$ The remainder of the proof is reduced to showing the bound
\begin{equation} \label{splus}
S_2:=\sum_{p^-< q < p^+} \sum_{\substack{f_1|D(p,q) \\ f_2|D(r,q) \\ (f_1f_2,2)=1}}\frac{L\left(1,\chi_1\right)
L\left(1,\chi_2\right)}{\varphi(f_1)\varphi(f_2)}\ll \frac{\sqrt{p}}{\log p}.
\end{equation}

Let $S_2'$ denote the double sum on the LHS of $\eqref{splus}$ with  $L\left(1,\chi_i;z^{8\alpha^2}\right)$
in place of $L\left(1,\chi_i\right)$ for $i=1,2$, where $z:=\log(4p)$ and $\alpha$ is a parameter $\geq 10$. We estimate the error term $S_2-S_2'$ by applying Proposition \ref{lemtwothree} once for $L\left(1,\chi_1\right)$ with $Q=4p$ and once for $L\left(1,\chi_2\right)$ with $Q=4r$. We have that $0\leq -D(p,q)\leq 4p$
 and $0\leq -D(r,q) \leq 4r$ for $q\in (p^-,p^+)$. Moreover, $\QQ(\sqrt{(2f_1)^2D(p,q)})=\QQ(\sqrt{D(p,q)}).$ If the conductor of $\chi_1$, which is the discriminant of $\QQ(\sqrt{D(p,q)}),$ 
 does not belong to the set ${\mathcal{E}}_\alpha(4p)$, or if the conductor of $\chi_2$, which is the discriminant of  $\QQ(\sqrt{D(r,q)}),$ 
 does not belong to ${\mathcal{E}}_\alpha(4r)$, we can bound $L\left(1,\chi_i\right)$ by $\log z$ from Mertens' theorem. For the exceptional sets ${\mathcal{E}}_\alpha(4p)$ and ${\mathcal{E}}_\alpha(4r)$ we use the convexity bound 
$L\left(1,\chi_i\right)\ll z$  for $i=1,2,$ respectively. This yields the estimate
\begin{align}
&S_2-S_2'\nonumber\\
\ll_\alpha& \frac{(\log z)^2}{z^\alpha}
\sum_{\substack{p^-< q < p^+ \\ {\rm disc}(\QQ(\sqrt{D(p,q)})) \not \in {\mathcal{E}}_\alpha(4p)
\\ {\rm disc}(\QQ(\sqrt{D(r,q)})) \not \in {\mathcal{E}}_\alpha(4r) }} 
\sum_{\substack{f_1|D(p,q)\\f_2|D(r,q) \\ (f_1f_2,2)=1}}\frac{1}{\varphi(f_1)\varphi(f_2)}\nonumber\\
+&z\log z\Bigg(\sum_{\substack{p^-< q < p^+ \\ 
{\rm disc}(\QQ(\sqrt{D(p,q)})) \in {\mathcal{E}}_\alpha(4p) \\ {\rm disc}(\QQ(\sqrt{D(r,q)})) \not 
\in {\mathcal{E}}_\alpha(4r)}} \sum_{\substack{f_1|D(p,q)\\f_2|D(r,q) \\ (f_1f_2,2)=1}}\frac{1}{\varphi(f_1)\varphi(f_2)}\nonumber\\
+&\sum_{\substack{p^-< q < p^+ \\ {\rm disc}(\QQ(\sqrt{D(p,q)})) \not \in {\mathcal{E}}_\alpha(4p) \\ {\rm disc}(\QQ(\sqrt{D(r,q)})) 
\in {\mathcal{E}}_\alpha(4r)}} \sum_{\substack{f_1|D(p,q)\\f_2|D(r,q) \\ (f_1f_2,2)=1}}\frac{1}{\varphi(f_1)\varphi(f_2)} \Bigg)+z^2\sum_{\substack{p^-< q < p^+ \\ {\rm disc}(\QQ(\sqrt{D(p,q)})) \in {\mathcal{E}}_\alpha(4p) \\ {\rm disc}(\QQ(\sqrt{D(r,q)})) \in 
{\mathcal{E}}_\alpha(4r)}} \sum_{\substack{f_1|D(p,q)\\f_2|D(r,q) \\ (f_1f_2,2)=1}}\frac{1}{\varphi(f_1)\varphi(f_2)}.\label{stwoerror}
\end{align}
For $q \in (p^-,p^+)$ such that $\Delta:= {\rm disc}(\QQ(\sqrt{D(p,q)})) \in {\mathcal{E}}_\alpha(4p)$ we have that $D(p,q)=\Delta m^2$ for 
some $m \in \NN$. Equivalently $(p+1-q)^2-\Delta m^2=4p,$ where $\Delta\equiv D(p,q)\equiv 1 \imod{4}$. Let $n=p+1-q$, 
then for a fixed $\Delta \in {\mathcal{E}}_\alpha(4p)$ we need to determine the quantity
\begin{align*}
r(4p,2):=&\#\{(m,n)\in \ZZ^2: n^2-\Delta m^2=4p \}, \\
=&\#\left\{\frac{n+m\sqrt{\Delta}}{2}\in {\mathcal{O}}_K : N\left(\frac{n+m\sqrt{\Delta}}{2}\right)=p\right\},
\end{align*}
where $K=\QQ(\sqrt{\Delta})=\QQ(\sqrt{D(p,q)}),{\mathcal{O}}_K $ is its ring of integers and $N(\cdot)$ is the norm of an element in $K$. 

Note that $$\#\{I\subseteq {\mathcal{O}}_K: N(I)=d\}=\left(1*\left(\frac{\Delta}{\cdot}\right)\right)(d),$$
where $N(I)$ denotes the norm of an ideal $I\subseteq {\mathcal{O}}_K$. Thus, $$\frac{r(4p,2)}{6}\leq \#\{I\subseteq {\mathcal{O}}_K: N(I)=p\}=\left(1*\left(\frac{\Delta}{\cdot}\right)\right)(p)$$
by the above equality. Hence, we conclude that $$r(4p,2)\leq6\sum_{k\mid p}\left(\frac{\Delta}{k}\right)\leq 12.$$
So there are at most 12 admissible pairs $(m,n)$ and therefore there are at most 12 admissible values of $q$ since $p$ is fixed. Thus,
$$\#\{p^-< q< p^+ : {\rm disc}(\QQ(\sqrt{D(p,q)})) \in {\mathcal{E}}_\alpha(4p)\} \leq 12\#{\mathcal{E}}_\alpha(4p)\ll p^\frac{1}{5},$$ 
since $\alpha\geq 10$. Similarly, we have that
$$\#\{p^-< q< p^+ : {\rm disc}(\QQ(\sqrt{D(r,q)})) \in {\mathcal{E}}_\alpha(4r)\} \leq 12\#{\mathcal{E}}_\alpha(4r)\ll r^\frac{1}{5} 
\ll p^\frac{1}{5}$$ and
\begin{align*}
&\#\Big\{p^-< q< p^+ : {\rm disc}(\QQ(\sqrt{D(p,q)})) \in {\mathcal{E}}_\alpha(4p)\:{\rm and}\:{\rm disc}(\QQ(\sqrt{D(r,q)})) \in 
{\mathcal{E}}_\alpha(4r)\Big\}\\
 &\leq 12\min\{\#{\mathcal{E}}_\alpha(4p),\#{\mathcal{E}}_\alpha(4r)\}\ll p^\frac{1}{5}.
\end{align*}

Since $f_1 \leq |D(p,q)|$ we have that $$\log\log f_1 \leq (\log\log |D(p,q)|) \ll 
\log \log p \ll \log z.$$ Thus, employing the bound $\displaystyle{\frac{1}{\varphi(f_1)}\ll \frac{\log \log f_1}{f_1}}$ yields
\begin{align}
\sum_{\substack{f_1|D(p,q) \\ (f_1,2)=1}}\frac{1}{\varphi(f_1)} &\ll \log z\sum_{\substack{f_1|D(p,q) \\ (f_1,2)=1}}\frac{1}{f_1} = 
\log z\prod_{\substack{ \ell\mid D(p,q) \\ \ell\neq 2}}\left(1-\frac{1}{\ell}\right)^{-1} \ll (\log z)^2.
\label{varphibounds}
\end{align}
The result is analogous for $D(r,q)$ and then applying the bounds on the exceptional set and the bound from
$\eqref{varphibounds}$ in $\eqref{stwoerror}$ yields
\begin{equation*}
S_2-S_2'\ll_\alpha \frac{\sqrt{p}(\log z)^6}{z^{1+\alpha}}+p^{\frac{1}{5}}(z(\log z)^5+z^2(\log z)^4),
\end{equation*}
and since $\alpha\geq 10$ we conclude that $S_2-S_2'\ll_\alpha \displaystyle{\frac{\sqrt{p}}{\log p}}$. 
Thus, it remains to show that
\begin{equation*}
S_2':=\sum_{p^-< q < p^+} \sum_{\substack{f_1|D(p,q) \\ f_2 |D(r,q)\\ (f_1f_2,2)=1}}\frac{L\left(1,\chi_1;z^{8\alpha^2}\right)
L\left(1,\chi_2;z^{8\alpha^2}\right)}{\varphi(f_1)\varphi(f_2)} \ll \frac{\sqrt{p}}{\log p}.
\end{equation*}
In order to do this, we find an upper bound for $L\left(1,\chi_1;z^{8\alpha^2}\right)$. Recall that
$$\chi_1:=\left(\frac{(2f_1)^2D(p,q)}{\cdot}\right).$$ By Mertens' theorem, we have that
\begin{align}
L(1,\chi_1;z^{8\alpha^2})&=\prod_{\ell\leq \sqrt{z}}\left(1-\frac{\chi_1(\ell)}{\ell}\right)^{-1}\prod_{\sqrt{z}\leq
 \ell\leq z^{8\alpha^2}}\left(1-\frac{\chi_1(\ell)}{\ell}\right)^{-1}\nonumber \\
&\ll_\alpha \prod_{\ell \leq \sqrt{z}}\left(1-\frac{\left(\frac{D(p,q)}{\ell}\right)}
{\ell}\right)^{-1}\prod_{\substack{\ell \leq \sqrt{z} 
\\ \ell \mid 2f_1}}\left(1-\frac{\left(\frac{D(p,q)}{\ell}\right)}{\ell}\right)\nonumber\\ 
 &\ll_\alpha  \frac{f_1}{\varphi(f_1)}
\prod_{\ell \leq \sqrt{z}}\left(1+\frac{\left(\frac{D(p,q)}{\ell}\right) }{\ell}\right),\label{lshortfour}
\end{align}
and similarly,
\begin{equation}\label{lshortfive}
L(1,\chi_2;z^{8\alpha^2}) \ll_\alpha \frac{f_2}{\varphi(f_2)} \prod_{\ell \leq \sqrt{z}}\left(1+
\frac{\left(\frac{D(r,q)}{\ell}\right) }{\ell}\right).
\end{equation}
Since the products on the RHS of $\eqref{lshortfour}$ and $\eqref{lshortfive}$ no longer depend on $f_1$ and $f_2$ we swap the sum and 
product to obtain the upper bound
\begin{align}
S_2'\ll& \sum_{p^-< q < p^+} \prod_{\ell \leq \sqrt{z}}\left(1+\frac{\left(\frac{D(p,q)}{\ell}\right)}
{\ell}\right)\left(1+\frac{\left(\frac{D(r,q)}{\ell}\right) }{\ell}\right) \sum_{\substack{f_1|D(p,q) 
\\ (f_1,2)=1}}\frac{f_1}{\varphi^2(f_1)}\sum_{\substack{f_2 |D(r,q)\\ 
(f_2,2)=1}}\frac{f_2}{\varphi^2(f_2)}.\label{splustwo}
\end{align}

We first consider the sum over $f_1$. Since $\displaystyle{\frac{f_1}{\varphi^2(f_1)}}$ is multiplicative by Mertens' theorem we have that
\begin{align}
\sum_{\substack{f_1|D(p,q) \\ (f_1,2)=1}}\frac{f_1}{\varphi^2(f_1)}&=\prod_{\substack{\ell |D(p,q) \\\ell\neq 2}}\left(1+\frac{\ell^2}{(\ell-1)^3}\right)
\ll\prod_{\substack{\ell |D(p,q) \\\ell\nmid 2f_2}}\left(1+\frac{1}{\ell}\right)\prod_{\substack{\ell |(D(p,q),f_2) \\ \ell\neq 2}}\left(1+\frac{1}{\ell}\right)\nonumber\\
&\ll \frac{f_2}{\varphi(f_2)}\prod_{\substack{\ell | D(p,q) \\ \ell \nmid 2f_2}}\left(1+\frac{1}{\ell}\right)
=\frac{f_2}{\varphi(f_2)}\prod_{\substack{\ell|D(p,q) \\ \ell \nmid 2f_2 \\ \ell \leq z^\alpha}}\left(1+\frac{1}{\ell}\right)(1+O(z^{-\alpha+1}))\nonumber \\
&\ll \frac{f_2}{\varphi(f_2)}\prod_{\substack{\ell|D(p,q) \\ \ell \nmid 2f_2 \\ \ell \leq \sqrt{z}}}\left(1+\frac{1}{\ell}\right)=\frac{f_2}{\varphi(f_2)}
\sum_{\substack{f_1 | D(p,q) \\ (f_1,2f_2)=1 \\ P^+(f_1)\leq \sqrt{z}}}\frac{\mu^2(f_1)}{f_1}. \label{coprimefone}
\end{align}
Replacing the RHS of $\eqref{coprimefone}$ in $\eqref{splustwo}$ yields
\begin{equation}\label{splustwoprime}
S_2'\ll \sum_{p^-< q < p^+} \prod_{\ell \leq \sqrt{z}}\left(1+\frac{\left(\frac{D(p,q)}{\ell}\right)}
{\ell}\right)\left(1+\frac{\left(\frac{D(r,q)}{\ell}\right) }{\ell}\right) 
\sum_{\substack{f_1 | D(p,q) \\ (f_1,2)=1 \\ P^+(f_1)\leq \sqrt{z}}}\frac{\mu^2(f_1)}{f_1}
\sum_{\substack{f_2 |D(r,q)\\ (f_2,2f_1)=1}}\frac{f_2^2}{\varphi^3(f_2)}.
\end{equation}

As in $\eqref{coprimefone}$ we have that
\begin{equation}
\sum_{\substack{f_2|D(r,q) \\ (f_2,2f_1)=1}}\frac{f_2^2}{\varphi^3(f_2)}
=\prod_{\substack{ \ell |D(r,q) \\ (\ell,2f_1)=1}}\left(1+\frac{\ell^3}{(\ell-1)^4}\right)
\ll\sum_{\substack{f_2 | D(r,q) \\ (f_2,2f_1)=1 \\ P^+(f_2)\leq \sqrt{z}}}\frac{\mu^2(f_2)}{f_2}. \label{coprimeftwo}
\end{equation}

Replacing $\eqref{coprimeftwo}$ in $\eqref{splustwoprime}$ yields
\begin{align}
S_2'&\ll \sum_{p^-< q < p^+} \prod_{\ell \leq \sqrt{z}}\left(1+\frac{\left(\frac{D(p,q)}{\ell}\right)}
{\ell}\right)\left(1+\frac{\left(\frac{D(r,q)}{\ell}\right) }{\ell}\right) 
\sum_{\substack{f_1 | D(p,q) \\ (f_1,2)=1 \\ P^+(f_1)\leq \sqrt{z}}}\frac{\mu^2(f_1)}{f_1}
\sum_{\substack{f_2 | D(r,q) \\ (f_2,2f_1)=1 \\ P^+(f_2)\leq \sqrt{z}}}\frac{\mu^2(f_2)}{f_2}. 
\label{splusthree}
\end{align}

Similar to $\eqref{lshortfour}$ we have that
\begin{align}
\prod_{\ell \leq \sqrt{z}}\left(1+\frac{\left(\frac{D(p,q)}{\ell}\right)}{\ell}\right)
\ll &\frac{f_1f_2}{\varphi(f_1)\varphi(f_2)}\prod_{\substack{\ell \leq \sqrt{z}\\\ell\nmid 2f_1f_2}}\left(1+\frac{\left(\frac{D(p,q)}{\ell}\right)}{\ell}\right), \label{lshortsix}
\end{align}
and
\begin{equation}
\prod_{\ell \leq \sqrt{z}}\left(1+\frac{\left(\frac{D(r,q)}{\ell}\right)}{\ell}\right)\ll \frac{f_1f_2}{\varphi(f_1)\varphi(f_2)}\prod_{\substack{\ell \leq \sqrt{z}\\\ell\nmid 2f_1f_2}}\left(1+\frac{\left(\frac{D(r,q)}{\ell}\right)}{\ell}\right). \label{lshortseven}
\end{equation}

Combining $\eqref{splusthree}, \eqref{lshortsix},$ and $\eqref{lshortseven}$ gives
\begin{align}
S_2'&\ll \sum_{\substack{P^+(f_1), P^+(f_2)\leq \sqrt{z} \\ (f_1,2)=(f_2,2f_1)=1 }}\frac{\mu^2(f_1)\mu^2(f_2)f_1f_2}{\varphi^2(f_1)\varphi^2(f_2)}\sum_{\substack{p^-< q < p^+ \\ f_1 | D(p,q) \\f_2 | D(r,q) }}
\prod_{\substack{\ell \leq \sqrt{z}\\\ell\nmid 2f_1f_2}}\left(1+\frac{\left(\frac{D(p,q)}{\ell}\right)}{\ell}\right)\left(1+\frac{\left(\frac{D(r,q)}{\ell}\right)}{\ell}\right). \label{splusfour}
\end{align}

We have that
\begin{equation}\label{eulerposione}
\prod_{\substack{\ell \leq \sqrt{z} \\ \ell \nmid 2f_1f_2}}\left(1+\frac{\left(\frac{D(p,q)}{\ell}\right)}{\ell}\right)= \sum_{\substack{P^+(n_1)\leq \sqrt{z} \\ (n_1,2f_1f_2)=1}}\frac{\mu^2(n_1)}{n_1}\left(\frac{D(p,q)}{n_1}\right),
\end{equation}
and likewise 
\begin{equation}\label{eulerpositwo}
\prod_{\substack{\ell \leq \sqrt{z} \\ \ell \nmid 2f_1f_2}}\left(1+\frac{\left(\frac{D(r,q)}{\ell}\right)}{\ell}\right)= \sum_{\substack{P^+(n_2)\leq \sqrt{z} \\ (n_2,2f_1f_2)=1}}\frac{\mu^2(n_2)}{n_2}\left(\frac{D(r,q)}{n_2}\right).
\end{equation}

Combining $\eqref{splusfour}$ with $\eqref{eulerposione}$ and $\eqref{eulerpositwo}$ and breaking up the RHS of $\eqref{splusfour}$ into sums over primes $q \mid 2f_1f_2n_1n_2$ and $q \nmid 2f_1f_2n_1n_2$ yields
\begin{align}
S_2'\ll&\sum_{\substack{P^+(f_1), P^+(f_2)\leq \sqrt{z} \\ (f_1,2)=(f_2,2f_1)=1 }}
 \frac{\mu^2(f_1)\mu^2(f_2)f_1f_2}{\varphi^2(f_1)\varphi^2(f_2)}\sum_{\substack{P^+(n_1),P^+(n_2) \leq \sqrt{z}\\(n_1n_2,2f_1f_2)=1}}
 \frac{\mu^2(n_1)\mu^2(n_2)}{n_1n_2}\nonumber \\
 \times&\sum_{\substack{p^-< q < p^+ \\ f_1 | D(p,q),f_2 | D(r,q) \\ (q,2f_1f_2n_1n_2)=1}}\left(\frac{D(p,q)}{n_1}\right)
\left(\frac{D(r,q)}{n_2}\right)\nonumber\\
+&\sum_{\substack{P^+(f_1), P^+(f_2)\leq \sqrt{z} \\ (f_1,2)=(f_2,2f_1)=1 }}
 \frac{\mu^2(f_1)\mu^2(f_2)f_1f_2}{\varphi^2(f_1)\varphi^2(f_2)}\sum_{\substack{P^+(n_1),P^+(n_2) \leq \sqrt{z}\\(n_1n_2,2f_1f_2)=1}} 
 \frac{\mu^2(n_1)\mu^2(n_2)}{n_1n_2}\nonumber \\
\times&\sum_{\substack{p^-< q < p^+ \\ f_1 | D(p,q),f_2 | D(r,q) \\ q\mid 2f_1f_2n_1n_2}}\left(\frac{D(p,q)}{n_1}\right)
\left(\frac{D(r,q)}{n_2}\right).\label{splusacop}
\end{align}

We have that the second sum in $\eqref{splusacop}$ is bounded by
\begin{equation}
\sum_{\substack{P^+(f_1), P^+(f_2)\leq \sqrt{z} \\ (f_1,2)=(f_2,2f_1)=1 }}
 \frac{\mu^2(f_1)\mu^2(f_2)f_1f_2\tau(f_1)\tau(f_2)}{\varphi^2(f_1)\varphi^2(f_2)}
\sum_{\substack{P^+(n_1),P^+(n_2) \leq \sqrt{z}\\(n_1n_2,2f_1f_2)=1}} \frac{\mu^2(n_1)\mu^2(n_2)\tau(n_1n_2)}{n_1n_2}, \label{sanotcop}
\end{equation}
where $\tau(n)$ denotes the number of divisors of $n$. We have that $\tau(n_1n_2)\leq \tau(n_1)\tau(n_2)$ and
\begin{align}
\sum_{\substack{P^+(n), \leq \sqrt{z}\\(n,2f_1f_2)=1}} \frac{\mu^2(n)\tau(n)}{n}&=\prod_{\substack{\ell \leq \sqrt{z} \\ \ell\nmid 2f_1f_2}}\left(1+\frac{2}{\ell}\right)
\ll (\log z)^2, \label{tauonebound}
\end{align}
by Mertens' theorem and similarly,
\begin{align}
\sum_{\substack{P^+(f), \leq \sqrt{z}\\(f,2)=1}} \frac{\mu^2(f)f\tau(f)}{\varphi^2(f)}&=\prod_{\substack{\ell \leq \sqrt{z} \\ \ell\nmid 2}}\left(1+\frac{2\ell}{(\ell-1)^2}\right)
\ll (\log z)^2. \label{tautwobound}
\end{align}

Thus, from $\eqref{tauonebound}$ and $\eqref{tautwobound}$ we have that $\eqref{sanotcop}$ is bounded by $(\log z)^8$ and we conclude that the second term in 
$\eqref{splusacop}$ is smaller than $\displaystyle{\frac{\sqrt{p}}{\log p}}$. Thus, it remains to show 
 \begin{align}
&\sum_{\substack{P^+(f_1), P^+(f_2)\leq \sqrt{z} \\ (f_1,2)=(f_2,2f_1)=1 }}
 \frac{\mu^2(f_1)\mu^2(f_2)f_1f_2}{\varphi^2(f_1)\varphi^2(f_2)}\sum_{\substack{P^+(n_1),P^+(n_2) \leq \sqrt{z}\\(n_1n_2,2f_1f_2)=1}}
 \frac{\mu^2(n_1)\mu^2(n_2)}{n_1n_2}\nonumber \\
 \times&\sum_{\substack{p^-< q < p^+ \\ f_1 | D(p,q),f_2 | D(r,q) \\ (q,2f_1f_2n_1n_2)=1}}\left(\frac{D(p,q)}{n_1}\right)
\left(\frac{D(r,q)}{n_2}\right)\ll \frac{\sqrt{p}}{\log p}. \label{splusfourcop}
\end{align}

Let $\lambda^+$ be the function defined in the fundamental lemma of sieve methods, Lemma \ref{lemtwoone} with $y=p^\frac{1}{6}$ and $D=y^2$. Then we have that the LHS of $\eqref{splusfourcop}$ is less than or equal to
\begin{align}
&\sum_{\substack{P^+(f_1), P^+(f_2)\leq \sqrt{z} \\ (f_1,2)=(f_2,2f_1)=1 }}
 \frac{\mu^2(f_1)\mu^2(f_2)f_1f_2}{\varphi^2(f_1)\varphi^2(f_2)}\sum_{\substack{P^+(n_1),P^+(n_2) \leq \sqrt{z}\\(n_1n_2,2f_1f_2)=1}}
 \frac{\mu^2(n_1)\mu^2(n_2)}{n_1n_2}\nonumber \\
 \times&\sum_{\substack{p^-\leq m\leq p^+ \\f_1\mid D(p,m),f_2 \mid D(r,m)\\ (m,2f_1f_2n_1n_2)=1}}(\lambda^+*1)(m)\left(\frac{D(p,m)}{n_1}\right)
 \left(\frac{D(r,m)}{n_2}\right),\label{sprimewithm}
\end{align}
by the positivity of the Euler product in $\eqref{eulerposione}$ and $\eqref{eulerpositwo}$. Hence, $\eqref{sprimewithm}$ becomes
\begin{align}
S_3:=&\sum_{\substack{P^+(f_1),P^+(f_2)\leq \sqrt{z} \\(f_1,2)=(f_2,2f_1)=1}}  \frac{\mu^2(f_1)\mu^2(f_2)f_1f_2}{\varphi^2(f_1)\varphi^2(f_2)}\sum_{\substack{P^+(n_1),P^+(n_2) \leq \sqrt{z}\\(n_1n_2,2f_1f_2)=1}} \frac{\mu^2(n_1)\mu^2(n_2)}{n_1n_2}\nonumber \\
\times&\sum_{\substack{a\leq D \\ (a,2f_1f_2n_1n_2)=1}}\lambda^+(a)\sum_{\substack{p^-< m < p^+ \\f_1\mid D(p,m), f_2 \mid D(r,m) \\ a\mid m}}\left(\frac{D(p,m)}{n_1}\right)\left(\frac{D(r,m)}{n_2}\right). \label{sprimefour}
\end{align}

Now we split the integers in the interval $m \in (p^-,p^+)$ according to the congruence class of $D(p,m)\imod{n_1}$ and $D(r,m)\imod{n_2}$. Thus, $\eqref{sprimefour}$ becomes
\begin{align*}
S_3=&\sum_{\substack{P^+(f_1),P^+(f_2)\leq \sqrt{z} \\(f_1,2)=(f_2,2f_1)=1}}  \frac{\mu^2(f_1)\mu^2(f_2)f_1f_2}{\varphi^2(f_1)\varphi^2(f_2)}\sum_{\substack{P^+(n_1),P^+(n_2) \leq \sqrt{z}\\(n_1n_2,2f_1f_2)=1}} \frac{\mu^2(n_1)\mu^2(n_2)}{n_1n_2}\nonumber \\
\times&\sum_{\substack{a\leq D \\ (a,2f_1f_2n_1n_2)=1}}\lambda^+(a)\sum_{\substack{b_1\in \ZZ/n_1\ZZ \\b_2\in \ZZ/n_2\ZZ} }\left(\frac{b_1}{n_1}\right)\left(\frac{b_2}{n_2}\right)S(a,f_1,f_2,n_1,n_2,b_1,b_2),
\end{align*}
where 
\begin{equation*}
S(a,f_1,f_2,n_1,n_2,b_1,b_2):= \#\left\{\begin{array}{ll}
 & D(p,m)\equiv 0 \imod{f_1}\\
& D(r,m)\equiv 0 \imod{f_2} \\ 
p^-< m < p^+ ;& D(p,m) \equiv b_1 \imod{n_1}  \\  &   D(r,m) \equiv b_2 \imod{n_2}\\  & m\equiv 0 \imod{a}\end{array}\right\}.
\end{equation*}
Since $a,f_1,f_2$,and $[n_1,n_2]$ are all coprime we have that
\begin{align}
S(a,f_1,f_2,n_1,n_2,b_1,b_2)&=\left(\frac{4\sqrt{p}}{af_1f_2[n_1,n_2]}\right)\#T(a,f_1,f_2,n_1,n_2,b_1,b_2)\nonumber\\
&+O(\#T(a,f_1,f_2,n_1,n_2,b_1,b_2)),\label{sprimeest}
\end{align}
where
\begin{equation*}
T(a,f_1,f_2,n_1,n_2,b_1,b_2):=\left\{\begin{array}{ll} & D(p,m)\equiv 0 \imod{f_1} \\ & D(r,m)\equiv 0 \imod{f_2} \\ m\in \ZZ/af_1f_2[n_1,n_2]\ZZ; & D(p,m) \equiv b_1 \imod{n_1}\\  & D(r,m) \equiv b_2 \imod{n_2} \\ & m\equiv 0 \imod{a}\end{array}\right\}.
\end{equation*}

Therefore, we have from $\eqref{sprimeest}$ that $\eqref{sprimefour}$ becomes
\begin{align}
S_3=&4\sqrt{p}\sum_{\substack{P^+(f_1),P^+(f_2)\leq \sqrt{z} \\(f_1,2)=(f_2,2f_1)=1}}  \frac{\mu^2(f_1)\mu^2(f_2)}{\varphi^2(f_1)\varphi^2(f_2)}\sum_{\substack{P^+(n_1),P^+(n_2) \leq \sqrt{z}\\(n_1n_2,2f_1f_2)=1}} \frac{\mu^2(n_1)\mu^2(n_2)}{n_1n_2[n_1,n_2]}\nonumber \\
\times&\sum_{\substack{a\leq D \\ (a,2f_1f_2n_1n_2)=1}}\frac{\lambda^+(a)}{a}\sum_{\substack{b_1\in \ZZ/n_1\ZZ \\b_2\in \ZZ/n_2\ZZ} }\left(\frac{b_1}{n_1}\right)\left(\frac{b_2}{n_2}\right)\#T(a,f_1,f_2,n_1,n_2,b_1,b_2)\nonumber\\
+&O\Bigg(\sum_{\substack{P^+(f_1),P^+(f_2)\leq \sqrt{z} \\(f_1,2)=(f_2,2f_1)=1}}  \frac{\mu^2(f_1)\mu^2(f_2)f_1f_2}{\varphi^2(f_1)\varphi^2(f_2)}\sum_{\substack{P^+(n_1),P^+(n_2) \leq \sqrt{z}\\(n_1n_2,2f_1f_2)=1}} \frac{\mu^2(n_1)\mu^2(n_2)}{n_1n_2}\nonumber \\
\times&\sum_{\substack{a\leq D \\ (a,2f_1f_2n_1n_2)=1}}|\lambda^+(a)|\sum_{\substack{b_1\in \ZZ/n_1\ZZ \\b_2\in \ZZ/n_2\ZZ} }\#T(a,f_1,f_2,n_1,n_2,b_1,b_2)\Bigg). \label{sprimefive}
\end{align}

By the Chinese remainder theorem we have that $$\#T(a,f_1,f_2,n_1,n_2,b_1,b_2)=\#T(a)\#T(f_1)\#T(f_2)\prod_{\ell \mid [n_1,n_2]}\#T^{(\ell)}(n_1,n_2,b_1,b_2),$$ where 
\begin{align}
\#T(a):=&\#\{m\in \ZZ/a\ZZ: m\equiv 0\imod{a}\}=1,\nonumber\\
\#T(f_1):=&\#\{m\in \ZZ/f_1\ZZ: D(p,m)\equiv 0\imod{f_1}\},\nonumber\\
\#T(f_2):=&\#\{m\in \ZZ/f_2\ZZ: D(r,m)\equiv 0\imod{f_2}\},\nonumber\\
\#T^{(\ell)}(n_1,n_2,b_1,b_2):=&\#\Bigg\{m\in \ZZ/\ell^{\nu_\ell([n_1,n_2])}\ZZ:  D(p,m)\equiv b_1\imod{\ell^{\nu_\ell(n_1)}} \nonumber\\
&\quad {\rm and}\: D(r,m)\equiv b_2\imod{\ell^{\nu_\ell(n_2)}}\Bigg\}. \label{tellsum}
 \end{align}

Note that $T(f_i)$ is multiplicative for $i=1,2$ and since we sum over odd, square-free $f_i$ in $\eqref{sprimefive}$ we have that
\begin{align}
\#T(f_1)=\prod_{\ell \mid f_1} \#\{m\in \ZZ/\ell\ZZ : (p+1-m)^2\equiv 4p\imod{\ell}\}=&\prod_{\ell \mid f_1}\left(1+\left(\frac{p}{\ell}\right)\right)= \sum_{d\mid f_1}\mu^2(d)\left(\frac{p}{d}\right)\label{huvbound}
\end{align}
 and similarly, $$\#T(f_2)=\prod_{\ell \mid f_2}\left(1+\left(\frac{r}{\ell}\right)\right)=\sum_{d\mid f_2}\mu^2(d)\left(\frac{r}{d}\right).$$

Thus, $\#T(f_i) \leq \tau(f_i)$ for all square-free integers $f_i$ for $i=1,2$. Now we consider the following function
\begin{equation*}
c(n_1,n_2):=\sum_{\substack{b_1\in \ZZ/n_1\ZZ \\ b_2 \in \ZZ/n_2 \ZZ}}\left(\frac{b_1}{n_1}\right)\left(\frac{b_2}{n_2}\right)\prod_{\ell \mid [n_1,n_2]}\#T^{(\ell)}(n_1,n_2,b_1,b_2).
\end{equation*}

Suppose that $n_1=n_1'n_1'', n_2=n_2'n_2''$ and $(n_1'n_2',n_1''n_2'')=1$. Then by the Chinese remainder theorem we have that 
\begin{align*}
 c(n_1'n_1'',n_2'n_2'')&=\sum_{\substack{b_1\in \ZZ/n_1'n_1''\ZZ \\ b_2\in \ZZ/n_2'n_2''\ZZ}}\left(\frac{b_1}{n_1'n_1''}\right)\left(\frac{b_2}{n_2'n_2''}\right)
 \prod_{\ell \mid [n_1'n_1'',n_2'n_2'']}\#T^{(\ell)}(n_1'n_1'',n_2'n_2'',b_1,b_2)\\
 &=\sum_{\substack{b_1'\in \ZZ/n_1'\ZZ \\ b_2'\in \ZZ/n_2'\ZZ}}\left(\frac{b_1'}{n_1'}\right)\left(\frac{b_2'}{n_2'}\right)
 \prod_{\ell \mid [n_1',n_2']}\#T^{(\ell)}(n_1',n_2',b_1',b_2')\\
 &\times\sum_{\substack{b_1''\in \ZZ/n_1''\ZZ \\ b_2''\in \ZZ/n_2''\ZZ}}\left(\frac{b_1''}{n_1''}\right)\left(\frac{b_2''}{n_2''}\right)
 \prod_{\ell \mid [n_1'',n_2'']}\#T^{(\ell)}(n_1'',n_2'',b_1'',b_2'')\\
 &=c(n_1',n_2')c(n_1'',n_2'').
\end{align*}
Thus, $c(n_1,n_2)$ is multiplicative and $[n_1'n_1'',n_2'n_2'']= [n_1',n_2'][n_1'',n_2'']$. We have that $n_1,n_2$ runs over square-free integers with $(n_1n_2,2f_1f_2)=1$ so it is enough to calculate $c(n_1,n_2)$ for primes $\ell \nmid 2f_1f_2$.
Since $c(1,1)=1$, we have three cases to consider, namely $c(\ell,1),c(1,\ell)$, and $c(\ell,\ell)$.

The cases $c(\ell,1)$ and $c(1,\ell)$ are completely similar and we have from $\eqref{tellsum}$ and $\eqref{huvbound}$ that
\begin{align*}
c(\ell,1)&=\sum_{b_1\in \ZZ/\ell\ZZ}\left(\frac{b_1}{\ell}\right)\#\{m\in \ZZ/\ell\ZZ : (p+1-m)^2\equiv 4p+b_1\imod{\ell}\}\\
&=\sum_{b_1\in \ZZ/\ell\ZZ}\left(\frac{b_1}{\ell}\right)\left(1+\left(\frac{4p+b_1}{\ell}\right)\right)= \sum_{b_1\in \ZZ/\ell\ZZ}\left(\frac{b_1}{\ell}\right)\left(\frac{4p+b_1}{\ell}\right)\\
&=\sum_{b_1\in \ZZ/\ell\ZZ}\left(\frac{b_1^2+4pb_1}{\ell}\right)=c(1,\ell).
\end{align*}
From \cite[Exercise 1.1.9]{SS:1} we have for $a\not \equiv 0 \imod{\ell}$ that
\begin{equation*}
\sum_{t \imod{\ell}}\left(\frac{at^2+bt+c}{\ell}\right)=\left\{
\begin{array}{ll} \left(\frac{a}{\ell}\right)(\ell-1) & {\rm
if}\: b^2-4ac \equiv 0 \imod{\ell}, \\
-\left(\frac{a}{\ell}\right) &{\rm if}\: b^2-4ac \not \equiv 0 \imod{\ell}.\end{array}\right.
\end{equation*}

Thus,
\begin{equation*}
c(\ell,1)=\left\{
\begin{array}{ll} \ell-1 & {\rm if}\: 16p^2 \equiv 0 \imod{\ell}, \\
-1 &{\rm if}\: 16p^2 \not \equiv 0 \imod{\ell}.\end{array}\right.
\end{equation*}
However, $\ell\nmid 2$ so if $16p^2 \equiv 0 \imod{\ell}$ then $\ell=p$. Since $P^+(n_1)\leq \sqrt{z}=\sqrt{\log 4p} <p$, 
we have that $c(\ell,1)=c(1,\ell)=-1$.

In the $c(\ell,\ell)$ case we have that
 \begin{equation*}
c(\ell,\ell)=\sum_{b_1,b_2\in \ZZ/\ell\ZZ}\left(\frac{b_1b_2}{\ell}\right)\#\{m\in \ZZ/\ell\ZZ: D(p,m)\equiv b_1\imod{\ell} \:{\rm and}\: D(r,m)\equiv b_2\imod{\ell}\}.
 \end{equation*}
We remark that there are at most two solutions to the equation $D(p,m)\equiv b_1\imod{\ell}$ since $D(p,m)$ is a 
quadratic polynomial in $m$. Let $m_0$ be one such solution. If $D(r,m_0)\not \equiv b_2\imod{\ell}$ 
then the two equations are not compatible. If $D(r,m_0) \equiv b_2\imod{\ell}$ then since the trace of 
$D(r,m)$ is fixed there will be at most $2$ values of $b_2$ that satisfy this equation. 
Hence, $$|c(\ell,\ell)|\leq\sum_{b_1\in \ZZ/\ell\ZZ}2 = 2\ell.$$

Combining the three cases, we conclude that 
$$|c(n_1,n_2)|\leq\prod_{\ell \mid (n_1,n_2)}|c(\ell,\ell)| \leq \prod_{\ell \mid (n_1,n_2)} 2\ell = 2^{\omega((n_1,n_2))}(n_1,n_2).$$
We now place our bounds from $\eqref{huvbound}$ and $c(n_1,n_2)$ into $\eqref{sprimefive}$ and we have that
\begin{align}
S_3\ll&\sqrt{p}\sum_{\substack{P^+(f_1),P^+(f_2)\leq \sqrt{z} \\(f_1,2)=(f_2,2f_1)=1}}  \frac{\mu^2(f_1)\mu^2(f_2)\tau(f_1)\tau(f_2)}{\varphi^2(f_1)\varphi^2(f_2)}\sum_{\substack{P^+(n_1),P^+(n_2) \leq \sqrt{z}\\(n_1n_2,2f_1f_2)=1}} \frac{\mu^2(n_1)\mu^2(n_2)(n_1,n_2)^2}{(n_1n_2)^{2-\epsilon}}\nonumber \\
\times&\Bigg|\sum_{\substack{a\leq D \\ (a,2f_1f_2n_1n_2)=1}}\frac{\lambda^+(a)}{a}\Bigg|+D\sum_{\substack{P^+(f_1),P^+(f_2)\leq \sqrt{z} \\(f_1,2)=(f_2,2f_1)=1}}  
\frac{\mu^2(f_1)\mu^2(f_2)f_1f_2\tau(f_1)\tau(f_2)}{\varphi^2(f_1)\varphi^2(f_2)}\nonumber\\
\times&\sum_{\substack{P^+(n_1),P^+(n_2) \leq \sqrt{z}\\(n_1n_2,2f_1f_2)=1}}
\frac{\mu^2(n_1)\mu^2(n_2)}{n_1n_2}\sum_{\substack{b_1\in \ZZ/n_1\ZZ \\b_2\in \ZZ/n_2\ZZ} }\prod_{\ell \mid [n_1,n_2]}
\#T^{(\ell)}(n_1,n_2,b_1,b_2). \label{sprimesix}
\end{align}
We first consider the second sum in $\eqref{sprimesix}$. Similarly to the function $c(n_1,n_2)$ defined above, the
function $$k(n_1,n_2):=\sum_{\substack{b_1\in \ZZ/n_1\ZZ \\b_2\in \ZZ/n_2\ZZ} }\prod_{\ell \mid [n_1,n_2]}
\#T^{(\ell)}(n_1,n_2,b_1,b_2)$$ is also multiplicative in $n_1$ and $n_2$. We have $k(1,1)=1,$ 
\begin{equation*}
k(\ell,1)=\sum_{b_1\in \ZZ/\ell\ZZ}\left(1+\left(\frac{4p+b_1}{\ell}\right)\right)=\ell=k(1,\ell),
\end{equation*}
and as in the case $c(\ell,\ell)$ above, we have that
$|k(\ell,\ell)|\leq\sum_{b_1\in \ZZ/\ell\ZZ}2 = 2\ell.$ Thus,

$$|k(n_1,n_2)|\leq\prod_{\ell \mid [n_1,n_2]}|k(\ell,1)k(1,\ell)k(\ell,\ell)| \leq \prod_{\ell \mid [n_1,n_2]} 2\ell^3 = 
2^{\omega([n_1,n_2])}[n_1,n_2]^3.$$ Substituting the bound above in $\eqref{sprimesix}$ we have that 
\begin{equation*}
\sum_{\substack{P^+(n_1),P^+(n_2) \leq \sqrt{z}\\(n_1n_2,2f_1f_2)=1}}
\frac{\mu^2(n_1)\mu^2(n_2)}{n_1n_2}2^{\omega([n_1,n_2])}[n_1,n_2]^3 \ll z^{3+\epsilon},
\end{equation*}
for $\epsilon>0$. Then by Mertens' theorem, for $i=1,2$ we have that
\begin{align*}
&\sum_{\substack{P^+(f_i)\leq \sqrt{z} \\(f_i,2)=1}} \frac{\mu^2(f_i)\tau(f_i)f_i}{\varphi^2(f_i)}\ll \frac{\sqrt{z}}{\log z}\prod_{\ell\leq \sqrt{z}}\left(1+\frac{2}{(\ell-1)^2}\right)\ll \frac{\sqrt{z}}{\log z},
\end{align*}
and thus, the second term in $\eqref{sprimesix}$ is bounded by $ Dz^{4+\epsilon}$. Then from Lemma \ref{lemtwoone} we have that $\eqref{sprimesix}$ becomes
\begin{align}
S_3\ll&\sqrt{p}\sum_{\substack{P^+(f_1),P^+(f_2)\leq \sqrt{z} \\(f_1,2)=(f_2,2f_1)=1}}  \frac{\mu^2(f_1)\mu^2(f_2)\tau(f_1)\tau(f_2)}{\varphi^2(f_1)\varphi^2(f_2)}\sum_{\substack{P^+(n_1),P^+(n_2) \leq \sqrt{z}\\(n_1n_2,2f_1f_2)=1}} \frac{\mu^2(n_1)\mu^2(n_2)(n_1,n_2)^2}{n_1^{2-\epsilon}n_2^{2-\epsilon}}\nonumber\\
\times&\prod_{\substack{\ell\leq y \\ \ell \nmid 2f_1f_2n_1n_2}} \left(1-\frac{1}{\ell}\right)+Dz^{4+\epsilon}.\label{sprimeseven}
\end{align}
By Mertens' theorem we have that
$$\prod_{\substack{\ell\leq y \\ \ell \nmid 2f_1f_2n_1n_2}} \left(1-\frac{1}{\ell}\right) \ll \frac{f_1f_2n_1n_2}{\varphi(f_1)\varphi(f_2)\varphi(n_1)\varphi(n_2)\log y},$$
and therefore the first term in the RHS of $\eqref{sprimeseven}$ is bounded by
\begin{align*}
\frac{\sqrt{p}}{\log y}\sum_{\substack{P^+(f_1),P^+(f_2)\leq \sqrt{z} \\(f_1,2)=(f_2,2f_1)=1}} \frac{\mu^2(f_1)\mu^2(f_2)\tau(f_1)\tau(f_2)f_1f_2}{\varphi^3(f_1)\varphi^3(f_2)}\sum_{\substack{P^+(n_1),P^+(n_2)\leq \sqrt{z}\\(n_1n_2,2f_1f_2)=1}}\frac{\mu^2(n_1)\mu^2(n_2)(n_1,n_2)^2}{\varphi(n_1)\varphi(n_2)n_1^{1-\epsilon}n_2^{1-\epsilon}}.
\end{align*}

We have that
\begin{align*}
&\sum_{\substack{P^+(n_1),P^+(n_2)\leq \sqrt{z}\\(n_1n_2,2f_1f_2)=1}}\frac{\mu^2(n_1)\mu^2(n_2)(n_1,n_2)^2}{\varphi(n_1)\varphi(n_2)n_1^{1-\epsilon}n_2^{1-\epsilon}}\\
 \ll &\sum_{\substack{P^+(d)\leq \sqrt{z} \\(d,2f_1f_2)=1}}\frac{\mu^2(d)}{d^{2-2\epsilon}}\sum_{\substack{P^+(m_1),P^+(m_2)\leq \frac{\sqrt{z}}{d}\\ n_1=dm_1, n_2=dm_2 \\ (d,m_1m_2)=1 \\ (m_1m_2,2f_1f_2)=1}} \frac{\mu^2(m_1)\mu^2(m_2)(\log\log dm_1)(\log\log dm_2)}{m_1^{2-\epsilon}m_2^{2-\epsilon}} \ll 1,
\end{align*}
and
\begin{align*}
&\sum_{\substack{P^+(f_1),P^+(f_2)\leq \sqrt{z} \\(f_1,2)=(f_2,2f_1)=1}} \frac{\mu^2(f_1)\mu^2(f_2)\tau(f_1)\tau(f_2)f_1f_2}{\varphi^3(f_1)\varphi^3(f_2)}\\
\ll & \sum_{\substack{P^+(f_1)\leq \sqrt{z} \\ (f_1,2)=1}} \frac{\mu^2(f_1)\tau(f_1)(\log \log f_1)^3}{f_1^2}\sum_{\substack{P^+(f_2)\leq \sqrt{z} \\(f_2,2f_1)=1}} \frac{\mu^2(f_2)\tau(f_2)(\log \log f_2)^3}{f_2^2}\ll 1.
\end{align*}

Thus, we conclude that
$$S_2\ll S_2'\ll S_3\ll\frac{\sqrt{p}}{\log y} + D(\log 4p)^{4+\epsilon} \ll \frac{\sqrt{p}}{\log p},$$
for $y=p^\frac{1}{6}, D=(p^\frac{1}{6})^2=p^{\frac{1}{3}},$ which completes the proof.
\end{proof}

The proof of Proposition \ref{hsinglebound} follows completely analogously to the steps taken in Proposition \ref{hsievebound}
and is essentially a special case of Chandee, David, Koukoulopoulos and Smith \cite[Proposition 4.1]{CDKS:1}.

\section{A short length of the average}
\begin{proof} (Proof of Lemma \ref{errorimprove})
Let $\chi_i$ and $\chi_i'$ be Dirichlet characters modulo $p_i$ for $1\leq i\leq L$ and let $\chi_0$ denote the principal character modulo $n$ for any integer $n$. For a Dirichlet character $\chi \imod{n}$, let $\bar{\chi}$ denote its complex conjugate of $\chi$ and let $${\mathcal{A}}(\chi):=\sum_{|a|\leq A}\chi(a) \quad {\rm and}\quad {\mathcal{B}}(\chi):=\sum_{|b|\leq B}\chi(b).$$ 
We recall from $\eqref{rpst}$ that  $R(P,S,T)$ is the number of integers $|a|\leq A,|b|\leq B$ such that there exists a vector
$(u_1,\ldots,u_L) \in \FF_{p_1}^*\times\cdots\times\FF_{p_L}^* $ satisfying
\begin{equation*}
 a\equiv s_iu_i^4\imod{p_i}, \quad b\equiv t_iu_i^6\imod{p_i} 
\quad{\rm for}\:1\leq i\leq L.
\end{equation*}
For $P:=(p_1,\ldots,p_L), S:=(s_1,\ldots,s_L), T:=(t_1,\ldots, t_L)$, and $U:=(u_1,\ldots, u_L)$ we have that
\begin{align}
R(P,S,T)=&\sum_{\substack{ |a|\leq A, |b|\leq B \\ \exists\; U\in\FF(P)^*\\ a\equiv s_iu_i^4\imod{p_i}, 
b\equiv t_iu_i^6\imod{p_i}\\ 1\leq i \leq L}} 1 \nonumber\\
=&\frac{1}{2^{L}}\sum_{\substack{|a|\leq A \\ |b|\leq B}}\sum_{ U\in\FF(P)^*}\prod_{i=1}^L 
\Bigg(\frac{1}{\varphi(p_i)^2}\sum_{\chi_i\imod{p_i}}\chi_i(s_iu_i^4)
\overline{\chi_i}(a)\sum_{\chi_i'\imod{p_i}}\chi_i'(t_iu_i^6)\overline{\chi_i'}(b)\Bigg)\nonumber\\
=&\frac{1}{2^{L}}\prod_{i=1}^L\frac{1}{(p_i-1)^2}\sum_{ U\in\FF(P)^*}\sum_{\substack{\chi_i, \chi_i'\imod{p_i} \\ 1\leq i \leq L}} \chi_i(s_i)\chi_i'(t_i)
\chi_i(u_i^4)\chi_i'(u_i^6)\nonumber\\
\times&\sum_{\substack{|a|\leq A \\ |b|\leq B}}\overline{\chi_1\cdots\chi_L}(a)\overline{\chi_1'\cdots\chi_L'}(b).\label{suv}
\end{align}
In $\eqref{suv}$ the factor $2^{-L}$ is present, since if there exists a $u_i\imod{p_i}$ such that $a\equiv s_iu_i^4\imod{p_i}$ and $b\equiv t_iu_i^6\imod{p_i}$ 
then there exists exactly two such $u_i$, namely $\pm u_i$.

By the orthogonality of Dirichlet characters, we have that the sum over $U$ becomes
\begin{equation}
\prod_{i=1}^L\sum_{ U\in\FF_{p_i}^*}\chi_i(u_i^4)\chi_i'(u_i^6)=\begin{cases}\prod_{i=1}^L (p_i-1) & \:{\rm if}\: \chi_i^4(\chi_i')^6=\chi_{0}\imod{p_i}\:{\rm for}\:1\leq i\leq L, \\ 0 & \:{\rm otherwise}.\end{cases} 
\label{suvtwo}
\end{equation}
Then from $\eqref{suv}$ and $\eqref{suvtwo}$ we have that

\begin{align}
R(P,S,T)=&\frac{1}{2^L}\sum_{\substack{\chi_1,\ldots,\chi_L \\ \chi_1',\ldots, \chi_L' \\ \chi_i^4(\chi_i')^6= \chi_{0} \imod{p_i} \\
{\rm for}\: 1\leq i\leq L}}\prod_{i=1}^L\left(\frac{\chi_i(s_i)\chi_i'(t_i)}{p_i-1}\right){\mathcal{A}}(\overline{\chi_1\cdots \chi_L})
{\mathcal{B}}(\overline{\chi_1'\cdots \chi_L'}) \nonumber\\
=&\frac{1}{2^L}\Bigg[\sum_{\substack{\chi_i=\chi_i'=\chi_0 \imod{p_i}\\
{\rm for}\: 1\leq i\leq L}} + \sum_{\substack{\chi_i=(\chi_i')^6=\chi_0 \imod{p_i}\\
{\rm for}\: 1\leq i\leq L\:{\rm and} \\ \exists 1\leq j \leq L \:{\rm s.t.}\: \chi_j'\neq \chi_0 \imod{p_j}}}+\sum_{\substack{\chi_i'=\chi_i^4=\chi_0 \imod{p_i}\\
{\rm for}\: 1\leq i\leq L\:{\rm and} \\ \exists 1\leq j \leq L \:{\rm s.t.}\: \chi_j\neq \chi_0 \imod{p_j}}}\nonumber \\
+&\sum_{\substack{\chi_i^4(\chi_i')^6=\chi_0 \imod{p_i}\\
{\rm for}\: 1\leq i\leq L\:{\rm and} \\ \exists 1\leq r,s \leq L \:{\rm s.t.}\: \chi_r\neq \chi_0 \imod{p_r},\\\chi_s'\neq \chi_0 \imod{p_s} }}\Bigg]\prod_{i=1}^L\left(\frac{\chi_i(s_i)\chi_i'(t_i)}{p_i-1}\right){\mathcal{A}}(\overline{\chi_1\cdots \chi_L})
{\mathcal{B}}(\overline{\chi_1'\cdots \chi_L'}).\label{thefourcs}
\end{align}
We denote the four sums in $\eqref{thefourcs}$ as follows,
\begin{align*}
R_1(P,S,T)&:=\frac{1}{2^L}\sum_{\substack{\chi_i=\chi_i'=\chi_0 \imod{p_i}\\
{\rm for}\: 1\leq i\leq L}}\prod_{i=1}^L\left(\frac{\chi_i(s_i)\chi_i'(t_i)}{p_i-1}\right){\mathcal{A}}(\overline{\chi_1\cdots \chi_L})
{\mathcal{B}}(\overline{\chi_1'\cdots \chi_L'}), \\
R_2(P,S,T)&:=\frac{1}{2^L}\sum_{\substack{\chi_i=(\chi_i')^6=\chi_0 \imod{p_i}\\
{\rm for}\: 1\leq i\leq L\:{\rm and} \\ \exists 1\leq j \leq L \:{\rm s.t.}\: \chi_j'\neq \chi_0 \imod{p_j}}}\prod_{i=1}^L\left(\frac{\chi_i(s_i)\chi_i'(t_i)}{p_i-1}\right){\mathcal{A}}(\overline{\chi_1\cdots \chi_L})
{\mathcal{B}}(\overline{\chi_1'\cdots \chi_L'}), \\
R_3(P,S,T)&:=\frac{1}{2^L}\sum_{\substack{\chi_i'=\chi_i^4=\chi_0 \imod{p_i}\\
{\rm for}\: 1\leq i\leq L\:{\rm and} \\ \exists 1\leq j \leq L \:{\rm s.t.}\: \chi_j\neq \chi_0 \imod{p_j}}}\prod_{i=1}^L\left(\frac{\chi_i(s_i)\chi_i'(t_i)}{p_i-1}\right){\mathcal{A}}(\overline{\chi_1\cdots \chi_L})
{\mathcal{B}}(\overline{\chi_1'\cdots \chi_L'}),\\
R_4(P,S,T)&:=\frac{1}{2^L}\sum_{\substack{\chi_i^4(\chi_i')^6=\chi_0 \imod{p_i}\\
{\rm for}\: 1\leq i\leq L\:{\rm and} \\ \exists 1\leq r,s \leq L \:{\rm s.t.}\: \chi_r\neq \chi_0 \imod{p_r},\\\chi_s'\neq \chi_0 \imod{p_s} }}\prod_{i=1}^L\left(\frac{\chi_i(s_i)\chi_i'(t_i)}{p_i-1}\right){\mathcal{A}}(\overline{\chi_1\cdots \chi_L})
{\mathcal{B}}(\overline{\chi_1'\cdots \chi_L'}).
\end{align*}

We recall the LHS of $\eqref{lemerrorbound}$, 
\begin{align*}
&\sumofpfracsing\sum_{\FFps}w(P,S,T) \left(R(P,S,T)-\frac{AB}{2^{L-2}p_1\cdots p_L}\right)\\
=&\sumofpfracsing\sum_{\FFps}w(P,S,T) \left(\sum_{j=1}^4R_j(P,S,T)-\frac{AB}{2^{L-2}p_1\cdots p_L}\right),
\end{align*}
by rewriting $R(P,S,T)$ as in $\eqref{thefourcs}$.

For $R_1(P,S,T)$ we have that $\chi_i=\chi_i'=\chi_0 \imod{p_i}$ for $1\leq i \leq L$ and hence,
\begin{align}
{\mathcal{A}}(\overline{\chi_1\cdots \chi_L})&=\sum_{|a|\leq A}\chi_0(a)=\sum_{\substack {|a|\leq A \\ (a,p_1\cdots p_L)=1}}1=2A\frac{\varphi(p_1\cdots p_L)}{p_1\cdots p_L}+O(\tau(p_1\cdots p_L))\nonumber \\
&=2A\left(\frac{(p_1-1)\cdots(p_L-1)}{p_1\cdots p_L}\right)+O_L(1)\label{achistuff}
\end{align}
and similarly, 
$${\mathcal{B}}(\overline{\chi_1'\cdots \chi_L'})=2B\left(\frac{(p_1-1)\cdots(p_L-1)}{p_1\cdots p_L}\right)+O_L(1).$$ 
Thus,
\begin{align} 
&R_1(P,S,T)\nonumber\\
=&\frac{1}{2^L}\prod_{j=1}^L \frac{1}{p_j-1}\left(\frac{2A(p_1-1)\cdots(p_L-1)}{p_1\cdots p_L}+O_L(1)\right)
\left(\frac{2B(p_1-1)\cdots(p_L-1)}{p_1\cdots p_L}+O_L(1)\right)\nonumber\\
=&\frac{AB}{2^{L-2}p_1\cdots p_L}+O_L\left(\frac{AB}{p^{L+\frac{1}{2}}}+\frac{A+B+1}{p^L}\right). \label{sjequalsone}
\end{align}

Recall from \eqref{dubsum} that
\begin{align}
\sum_{\FFps}w(P,S,T)=\prod_{i=1}^{L}(p_i-1)H(D(p_i,p_{i+1}))+O\left(p^{\frac{3L-1}{2}}(\log p)^{L-1}(\log \log p)^{L-1}\right).\label{dubsumtwo}
\end{align}

From $\eqref{sjequalsone}$ and  $\eqref{dubsumtwo}$ we have that 
\begin{align}
&\sumofpfracsing\sum_{\FFps}w(P,S,T) \left(R_1(P,S,T)-\frac{AB}{2^{L-2}p_1\cdots p_L}\right)\nonumber \\
\ll_L&\sumofpfracsing\left(\frac{AB}{p^{L+\frac{1}{2}}}+\frac{A+B+1}{p^L}\right)\prod_{j=1}^{L}(p_j-1)H(D(p_j,p_{j+1}))
\nonumber\\
\ll_L & \frac{AB (\log \log X)}{(\log X)^{L-1}}+\frac{(A+B+1)\sqrt{X}}{(\log X)^L},\label{jonebound}
\end{align}
by partial summation, Proposition \ref{hsievebound} and Proposition \ref{hsinglebound}. We have that $\eqref{jonebound}$ is smaller than the first two terms on the RHS in the error term of $\eqref{lemerrorbound}$. Thus, $\eqref{jonebound}$ is a lower order error term.

We now consider $R_2(P,S,T)$. From $\eqref{achistuff}$ we have that
\begin{align*}
R_2(P,S,T)&=\frac{1}{2^L} \sum_{\substack{(\chi_i')^6=\chi_0 \imod{p_i}\\
{\rm for}\: 1\leq i\leq L\:{\rm and} \\ \exists 1\leq j \leq L \:{\rm s.t.}\: \chi_j'\neq \chi_0 \imod{p_j}}} \prod_{j=1}^L\frac{\chi_j'(t_j)}{(p_j-1)}\left(2A\prod_{i=1}^L\frac{(p_i-1)}{p_i}+O_L(1)\right)
{\mathcal{B}}(\overline{\chi_1'\cdots \chi_L'})\\
&\ll_L\frac{A}{p_1\cdots p_L}\sum_{\substack{(\chi_i')^6=\chi_0 \imod{p_i}\\
{\rm for}\: 1\leq i\leq L\:{\rm and} \\ \exists 1\leq j \leq L \:{\rm s.t.}\: \chi_j'\neq \chi_0 \imod{p_j}}}|{\mathcal{B}}(\overline{\chi_1'\cdots \chi_L'})|.
\end{align*}

Similarly, we have that 
$$R_3(P,S,T)\ll_L\frac{B}{p_1\cdots p_L}\sum_{\substack{\chi_i^4=\chi_0 \imod{p_i}\\
{\rm for}\: 1\leq i\leq L\:{\rm and} \\ \exists 1\leq j \leq L \:{\rm s.t.}\: \chi_j\neq \chi_0 \imod{p_j}}}|{\mathcal{A}}(\overline{\chi_1\cdots \chi_L})|.$$ Thus, we have that 
\begin{align}
&\sumofpfracsing\sum_{\FFps}w(P,S,T)(R_2(P,S,T)+R_3(P,S,T))\nonumber \\
\ll_L&\sumofpfracsing\prod_{j=1}^{L}H(D(p_j,p_{j+1}))
\Bigg(A\sum_{\substack{(\chi_i')^6=\chi_0 \imod{p_i}\\
{\rm for}\: 1\leq i\leq L\:{\rm and} \\ \exists 1\leq j \leq L \:{\rm s.t.}\: \chi_j'\neq \chi_0 \imod{p_j}}} |{\mathcal{B}}(\overline{\chi_1'\cdots \chi_L'})|\nonumber\\
+&B\sum_{\substack{\chi_i^4=\chi_0 \imod{p_i}\\
{\rm for}\: 1\leq i\leq L\:{\rm and} \\ \exists 1\leq j \leq L \:{\rm s.t.}\: \chi_j\neq \chi_0 \imod{p_j}}}|{\mathcal{A}}(\overline{\chi_1\cdots \chi_L})|\Bigg).\label{errortermjtwo}
\end{align}

Let $$\displaystyle{\sumpchiprim}1=\sumofallps\sum_{\substack{(\chi_i')^6=\chi_0 \imod{p_i}\\
{\rm for}\: 1\leq i\leq L\:{\rm and} \\ \exists 1\leq j \leq L \:{\rm s.t.}\: \chi_j'\neq \chi_0 \imod{p_j}}}  1, $$ then by Holder's inequality we have that the first sum in $\eqref{errortermjtwo}$ becomes
\begin{align}
&A\sum_{\substack{p\leq X \\ p_i^-<p_{i+1}<p_i^+\\ 1\leq i \leq L-1}}\prod_{j=1}^L\frac{H(D(p_j,p_{j+1}))}{p_j}\sum_{\substack{(\chi_i')^6=\chi_0 \imod{p_i}\\
{\rm for}\: 1\leq i\leq L\:{\rm and} \\ \exists 1\leq j \leq L \:{\rm s.t.}\: \chi_j'\neq \chi_0 \imod{p_j}}}  |{\mathcal{B}}(\overline{\chi_1'\cdots \chi_L'})|\nonumber\\
\ll_L &A\Bigg(\sumpchiprim\prod_{j=1}^L\left(\frac{H(D(p_j,p_{j+1}))}{p_j}\right)^{\frac{2k}{2k-1}}\Bigg)^{1-\frac{1}{2k}}\Bigg(\sumpchiprim|{\mathcal{B}}(\overline{\chi_1'\cdots \chi_L'})|^{2k}\Bigg)^\frac{1}{2k}\label{sjholder}.
\end{align}

Since there are a bounded number of characters in the sums in $\eqref{sjholder}$ from $\eqref{hdellbound}$ we have that
\begin{align}
&\Bigg(\sumpchiprim\prod_{j=1}^L\left(\frac{H(D(p_j,p_{j+1}))}{p_j}\right)^{\frac{2k}{2k-1}}\Bigg)^{1-\frac{1}{2k}}\nonumber \\
\ll_L&\Bigg(\sum_{\substack{p\leq X \\ p_i^-<p_{i+1}<p_i^+\\ 1\leq i \leq L-1}}\left(\frac{(\log p)^L (\log \log p)^L}{p^\frac{L}{2}}\right)^{\frac{2k}{2k-1}}\Bigg)^{1-\frac{1}{2k}}
 \ll_L X^{\frac{1}{2}-\frac{L+1}{4k}}(\log X)^{\frac{L}{2k}}(\log\log X)^L.\label{holdboundone}
\end{align}

Let $J\subseteq\{1,\ldots,L\}$ be the set of positive integers such that if  $j\in J$ then $\chi_j'\neq \chi_0 \imod {p_j}$. For $R_2(P,S,T)$ we have that $J\neq \emptyset$. Thus, $$|{\mathcal{B}}(\overline{\chi_1'\cdots \chi_L'})|=\Bigg|\sum_{|b|\leq B}\overline{\chi_1'}(b)\cdots\overline{\chi_L'}(b) \Bigg|=\Bigg|\sum_{|b|\leq B}\prod_{j\in J}\overline{\chi_j'}(b)\prod_{j\not \in J}\overline{\chi_j'}(b)\Bigg|=\Bigg|\sum_{\substack{|b|\leq B \\ (b,\prod_{j\not \in J} p_j)=1}}\prod_{j\in J}\overline{\chi_j'}(b)\Bigg|.$$ Let $\tau_k(b;B)$ denote the number of ways of writing $b$ as a product of $k$ positive integers at most $B$. Then
$$\Bigg|\sum_{\substack{|b|\leq B \\ (b,\prod_{j\not \in J} p_j)=1}}\prod_{j\in J}\overline{\chi_j'}(b)\Bigg|^{2k}\ll_L\Bigg|\sum_{\substack{b\leq B \\ (b,\prod_{j\not \in J} p_j)=1}}\prod_{j\in J}\overline{\chi_j'}(b)\Bigg|^{2k}=\Bigg|\sum_{\substack{b\leq B^k \\ (b,\prod_{j\not\in J} p_j)=1}}\tau_k(b;B)\prod_{j\in J}\overline{\chi_j'}(b)\Bigg|^2.$$ 
Thus, for the second product in $\eqref{sjholder}$ we have that 
\begin{equation}
\Bigg(\sumpchiprim|{\mathcal{B}}(\overline{\chi_1'\cdots \chi_L'})|^{2k}\Bigg)^\frac{1}{2k}\ll_L \Bigg(\sumpchiprim\Bigg|\sum_{\substack{b\leq B^k \\ (b,\prod_{j\not\in J} p_j)=1}}\tau_k(b;B)\prod_{j\in J}\overline{\chi_j'}(b)\Bigg|^2\Bigg)^\frac{1}{2k}.\label{holdboundtwo}
\end{equation}

We have that $\prod_{j\in J}\overline{\chi_j'}(b)$ is a primitive character modulo $\prod_{j\in J} p_j$. Now we extend the sum in $\eqref{holdboundtwo}$ to a sum over all primitive characters modulo $d$ for all 
modulus $d\leq Q=X^{L}$, since $\prod_{j\in J} p_j \ll_L X^L$. 
Using the large sieve inequality, Theorem \ref{largesieve}, gives
\begin{align}
&\Bigg(\sumpchiprim\Bigg|\sum_{\substack{b\leq B^k \\ (b,\prod_{j\not\in J} p_j)=1}}\tau_k(b;B)\prod_{j\in J}\overline{\chi_j'}(b)\Bigg|^2\Bigg)^\frac{1}{2k}\nonumber\\
\ll_L &\Bigg(\sum_{\substack{d\leq X^L \\\chi\imod{d}\\ \chi \: {\rm primitive}}}\Bigg|\sum_{b\leq B^k}\tau_k(b;B)\chi(b)\Bigg|^2\Bigg)^\frac{1}{2k}
\ll_L \Bigg(\sum_{\substack{d\leq X^L \\\chi\imod{d}\\ \chi \: {\rm primitive}}}\Bigg|\sum_{b\leq B^k}\tau_k(b)\chi(b)\Bigg|^2\Bigg)^\frac{1}{2k}\nonumber\\
\ll_L &\Bigg((B^k+X^{2L})\sum_{b\leq B^k}\left|\tau_k(b)\right|^2\Bigg)^{\frac{1}{2k}} \ll_L \left((B^k+X^{2L})B^k\log^{k^2-1}(B^k)\right)^{\frac{1}{2k}}. \label{lalsieveboundprim}
\end{align}

Combining $\eqref{sjholder}, \eqref{holdboundone}$ and $\eqref{lalsieveboundprim}$ gives
\begin{align}
&A\sumofpfracsing\prod_{j=1}^{L}H(D(p_j,p_{j+1}))\sum_{\substack{(\chi_i')^6=\chi_0 \imod{p_i}\\
{\rm for}\: 1\leq i\leq L\:{\rm and} \\ \exists 1\leq j \leq L \:{\rm s.t.}\: \chi_j'\neq \chi_0 \imod{p_j}}}  |{\mathcal{B}}(\overline{\chi_1'\cdots \chi_L'})|\nonumber\\
\ll_L&A\left((B^k+X^{2L})B^k\log^{k^2-1}(B^k)\right)^{\frac{1}{2k}}X^{\frac{1}{2}-\frac{L+1}{4k}}(\log X)^{\frac{L}{2k}}(\log\log X)^L.
\label{lholdboundthreeprim}
\end{align}

First suppose that $B^k>X^{2L}$. Then we have that the RHS of $\eqref{lalsieveboundprim}$ becomes
\begin{align}
\left((B^k+X^{2L})B^k\log^{k^2-1}(B^k)\right)^{\frac{1}{2k}}
\ll_{k,L} B\log^{\frac{k^2-1}{2k}}B,\label{lholdboundfourprim}
\end{align}
for $k\geq 1$. Now suppose that $B^k\leq X^{2L}$ for all $k\geq 1$. Then we can replace $\log B$ by $\log X$ in $\eqref{lalsieveboundprim}$, which gives
\begin{align}
\left((B^k+X^{2L})B^k\log^{k^2-1}(B^k)\right)^{\frac{1}{2k}}
\ll_{k,L} \sqrt{B}X^\frac{L}{k}(\log X)^{\frac{k^2-1}{2k}}.\label{lholdboundfiveprim}
\end{align}

Since $$(B^k+X^{2L})^{\frac{1}{2k}}\ll_{k,L} \sqrt{B}+X^{\frac{L}{k}},$$ combining $\eqref{lholdboundfourprim}$ and $\eqref{lholdboundfiveprim}$ with $\eqref{lholdboundthreeprim}$ gives
\begin{align}
&A\sumofpfracsing\prod_{i=1}^{L}H(D(p_i,p_{i+1}))\sum_{\substack{(\chi_i')^6=\chi_0 \imod{p_i}\\
{\rm for}\: 1\leq i\leq L\:{\rm and} \\ \exists 1\leq j \leq L \:{\rm s.t.}\: \chi_j'\neq \chi_0 \imod{p_j}}} |{\mathcal{B}}(\overline{\chi_1'\cdots \chi_L'})|\nonumber \\
\ll_{L,k}&ABX^{\frac{1}{2}-\frac{L+1}{4k}}(\log X)^{\frac{L}{2k}}(\log\log X)^L\log^{\frac{k^2-1}{2k}}B+A\sqrt{B}X^{\frac{1}{2}+\frac{3L-1}{4k}}(\log X)^{\frac{k^2+L-1}{2k}}(\log\log X)^L \label{secondtermbound}.
\end{align}

Similarly we have that
\begin{align}
&B\sumofpfracsing\prod_{j=1}^{L}H(D(p_j,p_{j+1}))\sum_{\substack{\chi_i^4=\chi_0 \imod{p_i}\\
{\rm for}\: 1\leq i\leq L\:{\rm and} \\ \exists 1\leq j \leq L \:{\rm s.t.}\: \chi_j\neq \chi_0 \imod{p_j}}} |{\mathcal{A}}(\overline{\chi_1\cdots \chi_L})|\nonumber \\
\ll_{L,k}& ABX^{\frac{1}{2}-\frac{L+1}{4k}}(\log X)^{\frac{L}{2k}}(\log\log X)^L\log^{\frac{k^2-1}{2k}}A+B\sqrt{A}X^{\frac{1}{2}+\frac{3L-1}{4k}}(\log X)^{\frac{k^2+L-1}{2k}}(\log\log X)^L.\label{sjthreebound}
\end{align}

Thus, from $\eqref{secondtermbound}$ and $\eqref{sjthreebound}$ we have that $\eqref{errortermjtwo}$ becomes
\begin{align}
&\sumofpfracsing\sum_{\FFps}w(P,S,T)(R_2(P,S,T)+R_3(P,S,T))\nonumber \\
\ll_{k,L}& ABX^{\frac{1}{2}-\frac{L+1}{4k}}(\log X)^{\frac{L}{2k}}(\log\log X)^L(\log^{\frac{k^2-1}{2k}}A+\log^{\frac{k^2-1}{2k}}B)\nonumber\\
+&(A\sqrt{B}+B\sqrt{A})X^{\frac{1}{2}+\frac{3L-1}{4k}}(\log X)^{\frac{k^2+L-1}{2k}}(\log\log X)^L. \label{jtwosevenbound}
\end{align}

Now consider the final case $R_4(P,S,T)$. Let $$W(P,\chi_i,\chi_i'):= \sum_{\substack{1\leq s_i,t_i<p_i \\ 1\leq i\leq L}}w(P,S,T)\chi_i(s_i)\chi_i'(t_i).$$ Then we have that
\begin{align}
&\sumofpfracsing\sum_{\FFps}w(P,S,T)R_4(P,S,T)\nonumber\\
=&\frac{1}{2^L}\sumofallps \prod_{j=1}^L \frac{1}{p_j(p_j-1)}\sum_{\substack{\chi_i^4(\chi_i')^6=\chi_0 \imod{p_i}\\
{\rm for}\: 1\leq i\leq L\:{\rm and} \\ \exists 1\leq r,s \leq L \:{\rm s.t.}\: \chi_r\neq \chi_0 \imod{p_r},\\\chi_s'\neq \chi_0 \imod{p_s} }}W(P,\chi_i,\chi_i'){\mathcal{A}}(\overline{\chi_1\cdots \chi_L}){\mathcal{B}}(\overline{\chi_1'\cdots \chi_L'}).\label{sfoursumblah}
\end{align}
We use H{\"o}lder's inequality to obtain
\begin{align}
&\Bigg|\sum_{\substack{\chi_i^4(\chi_i')^6=\chi_0 \imod{p_i}\\
{\rm for}\: 1\leq i\leq L\:{\rm and} \\ \exists 1\leq r,s \leq L \:{\rm s.t.}\: \chi_r\neq \chi_0 \imod{p_r},\\\chi_s'\neq \chi_0 \imod{p_s} }}W(P,\chi_i,\chi_i'){\mathcal{A}}(\overline{\chi_1\cdots \chi_L}){\mathcal{B}}(\overline{\chi_1'\cdots \chi_L'})\Bigg| \nonumber\\
\leq&\Bigg|\sum_{\substack{\chi_i^4(\chi_i')^6=\chi_0 \imod{p_i}\\
{\rm for}\: 1\leq i\leq L\:{\rm and} \\ \exists 1\leq r,s \leq L \:{\rm s.t.}\: \chi_r\neq \chi_0 \imod{p_r},\\\chi_s'\neq \chi_0 \imod{p_s} }}\left|W(P,\chi_i,\chi_i')\right|^2\Bigg|^\frac{1}{2}\Bigg(\sum_{\substack{\chi_i^4(\chi_i')^6=\chi_0 \imod{p_i}\\
{\rm for}\: 1\leq i\leq L\:{\rm and} \\ \exists 1\leq r,s \leq L \:{\rm s.t.}\: \chi_r\neq \chi_0 \imod{p_r},\\\chi_s'\neq \chi_0 \imod{p_s} }}\left|{\mathcal{A}}(\overline{\chi_1\cdots \chi_L})\right|^4\Bigg)^{\frac{1}{4}} \nonumber\\
\times&\Bigg(\sum_{\substack{\chi_i^4(\chi_i')^6=\chi_0 \imod{p_i}\\
{\rm for}\: 1\leq i\leq L\:{\rm and} \\ \exists 1\leq r,s \leq L \:{\rm s.t.}\: \chi_r\neq \chi_0 \imod{p_r},\\\chi_s'\neq \chi_0 \imod{p_s} }}\left|{\mathcal{B}}(\overline{\chi_1'\cdots \chi_L'})\right|^4\Bigg)^\frac{1}{4}.\label{jfourcauchy}
\end{align}

We can extend the sums in the last two products in $\eqref{jfourcauchy}$ to a sum over all non-principal characters modulo $p_1\cdots p_L$. Thus, from Theorem \ref{fourthpower} we have that
\begin{align}
&\Bigg(\sum_{\substack{\chi_i^4(\chi_i')^6=\chi_0 \imod{p_i}\\
{\rm for}\: 1\leq i\leq L\:{\rm and} \\ \exists 1\leq r,s \leq L \:{\rm s.t.}\: \chi_r\neq \chi_0 \imod{p_r},\\\chi_s'\neq \chi_0 \imod{p_s} }}\left|{\mathcal{A}}(\overline{\chi_1\cdots \chi_L})\right|^4\sum_{\substack{\chi_i^4(\chi_i')^6=\chi_0 \imod{p_i}\\
{\rm for}\: 1\leq i\leq L\:{\rm and} \\ \exists 1\leq r,s \leq L \:{\rm s.t.}\: \chi_r\neq \chi_0 \imod{p_r},\\\chi_s'\neq \chi_0 \imod{p_s} }}\left|{\mathcal{B}}(\overline{\chi_1'\cdots \chi_L'})\right|^4\Bigg)^\frac{1}{4} \nonumber \\
\ll_L&\Bigg(\sum_{\chi\neq \chi_0 \imod{p_1\cdots p_L}}\Bigg|\sum_{|a|\leq A}\overline{\chi}(a)\Bigg|^4\Bigg)^{\frac{1}{4}}\Bigg(\sum_{\chi'\neq \chi_0 \imod{p_1\cdots p_L}}\Bigg|\sum_{|b|\leq B}\overline{\chi'}(b)\Bigg|^4\Bigg)^\frac{1}{4}\nonumber \\
\ll_L &\sqrt{ABp_1\cdots p_L}(\log p_1\cdots p_L)^3 \ll_L  \sqrt{ABp_1\cdots p_L}(\log p)^3.\label{jfourone}
\end{align}

Set $S':=(s_1',\ldots,s_L')$ and $T':=(t_1',\ldots,t_L').$
We then extend the first sum in $\eqref{jfourcauchy}$ to a sum over all possible products of characters modulo $p_1\cdots p_L$ (including the trivial character). Then we use the bound from $\eqref{dubsumtwo}$ to obtain
\begin{align}
&\sum_{\substack{\chi_i^4(\chi_i')^6=\chi_0 \imod{p_i}\\
{\rm for}\: 1\leq i\leq L\:{\rm and} \\ \exists 1\leq r,s \leq L \:{\rm s.t.}\: \chi_r\neq \chi_0 \imod{p_r},\\\chi_s'\neq \chi_0 \imod{p_s} }}\left|W(P,\chi_i,\chi_i')\right|^2\leq \sum_{\substack {\chi_i,\chi_i'\imod{p_i}\\ 1\leq i \leq L}}\left|W(P,\chi_i,\chi_i')\right|^2\nonumber \\
&\leq \sum_{\FFps}\sum_{S',T'\in\FF(P)^*}
w(P,S,T)\overline{w(P,S',T')}\sum_{\chi_i\imod{p_i}}\chi_i(s_i)\overline{\chi_i}(s_i')\sum_{\chi_i'\imod{p_i}}\chi_i'(t_i)\overline{\chi_i'}(t_i')\nonumber\\
&=\prod_{i=1}^L (p_i-1)^2\sum_{\FFps}\left|w(P,S,T)\right|^2\nonumber\\
&=p^{3L}\prod_{i=1}^{L}H(D(p_i,p_{i+1}))+O_L\left(p^{\frac{7L-1}{2}}(\log p)^L(\log \log p)^{L}\right),\label{jfourtwo}
\end{align}
since $|w(P,S,T)|^2=w(P,S,T)$. 

By combining $\eqref{jfourcauchy}, \eqref{jfourone}$ and $\eqref{jfourtwo}$ we have that
 \begin{align}
&\Bigg|\sum_{\substack{\chi_i^4(\chi_i')^6=\chi_0 \imod{p_i}\\
{\rm for}\: 1\leq i\leq L\:{\rm and} \\ \exists 1\leq r,s \leq L \:{\rm s.t.}\: \chi_r\neq \chi_0 \imod{p_r},\\\chi_s'\neq \chi_0 \imod{p_s} }}W(P,\chi_i,\chi_i'){\mathcal{A}}(\overline{\chi_1\cdots \chi_L}){\mathcal{B}}(\overline{\chi_1'\cdots \chi_L'})\Bigg| \nonumber\\
 \ll_L &\sqrt{AB}p^{2L}(\log p)^3\prod_{i=1}^{L}(H(D(p_i,p_{i+1})))^2.\label{wpfourstuff}
 \end{align}
 
 Then substituting $\eqref{wpfourstuff}$ into $\eqref{sfoursumblah}$ gives
\begin{align}
&\sumofpfracsing\sum_{\FFps}w(P,S,T)R_4(P,S,T)\nonumber\\
\ll_L&\sqrt{AB}\sum_{p\leq X}(\log p)^3 \sum_{\substack{p_i^-<p_{i+1}<p_i^+ \\ 1\leq i \leq L-1}}\prod_{j=1}^L \sqrt{H(D(p_j,p_{j+1}))}.\label{jfourprecauchy}
\end{align}

To obtain a better error term, instead of using the bound from $\eqref{hdellbound}$ for $H(D(p_j,p_{j+1}))$, we use Cauchy-Schwarz, Proposition \ref{hsievebound} and Proposition \ref{hsinglebound} to bound the inner sum in $\eqref{jfourprecauchy}$. This yields
\begin{align}
&\sum_{\substack{p_i^-<p_{i+1}<p_i^+ \\ 1\leq i \leq L-1}}\prod_{j=1}^L \sqrt{H(D(p_j,p_{j+1}))}\nonumber \\
\ll_L&\Bigg(\prod_{i=1}^{L-2}\frac{\sqrt{p_i}}{\log p_i}\sum_{p_i^-< p_{i+1} < p_i^+}H(D(p_i,p_{i+1}))
\frac{\sqrt{p}}{\log p}\sum_{p_{L-1}^-< p_{L} < p_{L-1}^+}H(D(p_{L-1},p_L))H(D(p_L,p))\Bigg)^{\frac{1}{2}} \nonumber \\
\ll_L &\prod_{i=1}^{L-2}\left(\frac{p_i}{\log p_i}\cdot \frac{\sqrt{p_i}}{\log p_i}\right)^\frac{1}{2}\left(\frac{\sqrt{p}}{\log p}\cdot\frac{ p^\frac{3}{2}}{\log p}\right)^\frac{1}{2}
\ll_L \frac{p^\frac{3L-2}{4}}{(\log p)^{L-1}}. \label{cauchyhbound}
\end{align}

From $\eqref{jfourprecauchy}$ and $\eqref{cauchyhbound}$ we have that
\begin{align}
&\sumofpfracsing\sum_{\FFps}w(P,S,T)R_4(P,S,T)
 \ll_L\sqrt{AB}X^{\frac{3L+2}{4}}(\log X)^{3-L}.
 \label{jeightbound}
  \end{align}
 Combining $\eqref{jtwosevenbound}$ and $\eqref{jeightbound}$ gives the result.
\end{proof}

\end{document}